\documentclass[leqno,a4paper, 11pt]{article}

\usepackage{amsfonts, amssymb}

\font\eulersm=eusm10 at 11pt
\def\esm#1{\hbox{\eulersm {#1}}}

\def\C{{\cal C}}
\def\Q{{\bf Q}}
\def\T{{\esm{T}}}
\def\O{{\cal O}}

\def\M{{\esm{M}}}
\def\F{{\bf F}}
\def\P{{\bf P}}
\def\Z{{\bf Z}}
\def\H{{\rm H}}
\def\m{\mathfrak{m}}

\def\bigdownarrow{\vphantom{\bigg|}\Big\downarrow}
\newenvironment{condlist}{\vspace{1ex}
     \begin{list}{}{%
          \setlength{\topsep}{0pt}%
          \setlength{\parsep}{0pt}%
          \setlength{\itemsep}{0pt}}%
     }{\end{list}\vspace{1ex}}
\def\underrel#1#2{\mathrel{\mathop{#1}\limits_{#2}}}
\def\Spf{\mathrm{Spf}}
\def\Proj{\mathrm{Proj}}
\def\Sym{\mathrm{Sym}}
\def\Spec{\mathrm{Spec}}
\def\Aut{\mathrm{Aut}}
\def\Hom{\mathrm{Hom}}

\def\PGL{\mathrm{PGL}}

\def\Ends{\mathrm{Ends}}
\def\Ad{\mathrm{Ad}}
\def\til#1{\widetilde{#1}}% tilde
\def\ovl#1{\overline{#1}}% overline
\def\pf{{\indent\textit{Proof.}\ }}
\def\qed{\hfill$\square$}
\newcounter{para}[section]
\renewcommand{\thepara}{\thesection.\arabic{para}}
\renewcommand{\thesection}{\arabic{section}}
\renewcommand{\paragraph}{\refstepcounter{para}
\indent{\bf{\thepara}}}
\newcommand{\sectioning}{\refstepcounter{section}
\indent{\bf \thesection.}}

\newcommand{\sd}{\rtimes}
\newcommand{\an}{{\mbox{{\rm {\tiny an}}}}}
\newcommand{\A}{\mbox{\rm{Aut}}}

\newcommand{\vv}{\vspace{1ex}}

\begin{document}

{\footnotesize \noindent Running head: Equivariant deformation in positive characteristic \hfill August 24, 2001

\noindent Math.\ Subj.\ Class.\ (2000): 14G22, 14D15, 14H37}

\vspace{3cm}

\begin{center} 

{\Large Equivariant deformation of Mumford curves and 

\vv

of ordinary curves in positive characteristic}

\vv

{\sl by} Gunther Cornelissen {\sl and} Fumiharu Kato

\vv

\end{center}

\vv

\vv

\begin{quote} {\small {\bf Abstract.} We compute the dimension of the tangent space to, and the Krull dimension of the pro-representable hull of two deformation functors. The first one is the ``algebraic'' deformation functor of an ordinary curve $X$ over a field of positive characteristic with prescribed action of a finite group $G$, and the data are computed in terms of the ramification behaviour of $X \rightarrow G \backslash X$. The second one is the ``analytic'' deformation functor of a fixed embedding of a finitely generated discrete group $N$
in $PGL(2,K)$ over a non-archimedean valued field $K$, and the data are
computed in terms of the Bass-Serre representation of $N$ via a graph of groups. Finally, if $\Gamma$ is a free subgroup of $N$ such that $N$ is contained
in the normalizer of $\Gamma$ in $PGL(2,K)$, then the Mumford curve associated to $\Gamma$ becomes
equipped with an action of $N/\Gamma$, and we show that the algebraic 
functor deforming the latter action coincides with the analytic functor deforming the embedding of $N$.}
\end{quote}

\vv

\vv

{\bf Introduction.}

\vspace{1ex}
Equivariant deformation theory is the correct framework for formulating
and answering questions such as ``given a curve $X$ of genus $g$ 
over a field $k$ and a finite group of automorphisms $ \rho : G \hookrightarrow
\A(X) $ of $X$, in how
many ways can $X$ be deformed into another curve of the same genus on
which the same group of automorphisms still acts?'' The precise
meaning of this question (at least infinitesimally) is related to the deformation functor $D_{X,\rho}$ of the pair
$(X,\rho)$, which associates to any element $A$ of the category $\C_k$ of local Artinian $k$-algebras with residue field $k$
the set of isomorphism classes of liftings $(X^\sim,\rho^\sim,\phi^\sim)$, 
where $X^\sim$ is a smooth scheme of finite type over $A$, $\phi^\sim$ 
is an isomorphism of $X^\sim \otimes k$ with $X$, and where $\rho^\sim 
: G \rightarrow \A_A(X)$ lifts $\rho$ via $\phi^\sim$. In general, 
$D_{X,\rho}$ has a pro-representable hull $H_{X,\rho}$ in the sense
of Schlessinger (\cite{Sc68}). This means that there is a smooth 
map of functors $\Hom(H_{X,\rho},-) \rightarrow D_{X,\rho}$ that
induces an isomorphism on the level of tangent spaces, where
$H_{X,\rho}$ is a Noetherian complete local $k$-algebra
with residue field $k$.
The
above question can then be reformulated as the computation of the Krull-dimension of $H_{X,\rho}$.
 Two remarks are in order: if $g \geq 2$, then the functor is even 
pro-representable by $H_{X,\rho}$, so $\mathrm{Spf} \ H_{X,\rho}$ is a formal scheme which is, so 
to speak, the universal basis of a family of curves which 
have the same automorphism group as $X$. Secondly, $H_{X,\rho}$ is 
even algebrizable (cf. Grothendieck \cite{Gro}, \S 3), and the underlying
algebraic scheme over $k$ might be considered
as the genuine ``universal basis'' scheme.

If the characteristic of the ground field $k$ is zero, the dimension 
of $H_{X,\rho}$ is easy to compute. All obstructions and group cohomology (cf.\ infra) 
disappear, and we find that 

\begin{equation} \label{char0}
\dim H_{X,\rho} = 3 g_Y - 3 +n, \end{equation}

\noindent where $Y:=G \backslash X$ is the quotient of $X$, $g_Y$ its genus and $n$ is the number
of branch points on $Y$ (note that $3g_Y-3$ is the degree of freedom of varying 
the moduli of $Y$, and one extra degree of freedom comes in for every 
branch point). This result can be found in any classical 
text on Riemann surfaces ({\sl e.g.,} \cite{Farkas:80}, V.2.2); the moral
is that ramification data of $X \rightarrow Y$ provide all the necessary
information for computing $H_{X,\rho}$. 

In this work, we are interested in the corresponding question 
in positive characteristic. Let us first present the motivating example for our studies: moduli schemes for rank two Drinfeld modules with principal level structure (see, e.g., Gekeler \& Reversat \cite{Gekeler:96}).

\medskip

{\bf Example.} \ Let $q=p^t$, $F={\bf F}_q(T)$, and $A={\bf F}_q[T]$; let 
$F_\infty={\bf F}_q((T^{-1}))$  be the completion of $F$ and $C$ a completion of the algebraic closure of $F_\infty$. On Drinfeld's ``upper half plane'' $\Omega:= {\bf P}^1_C - {\bf P}^1_{F_\infty}$ (which is a rigid analytic space over $C$), the group $GL(2,A)$ acts by fractional transformations. Let ${\cal Z} \cong {\bf F}_q^*$ be its center. For ${\mathfrak n} \in A$, the quotients of $\Omega$ by congruence subgroups $\Gamma({\mathfrak n}) = \{ \gamma \in GL(2,A) : \gamma = {\bf 1} \mbox{ mod } {\mathfrak n} \}$ are open analytic curves which can be
compactified to projective curves $X({\mathfrak n})$ by adding finitely many cusps.
These curves are analogues in the function field setting of classical
modular curves $X(n)$ for $n \in {\bf Z}$. 
Clearly, elements from $G({\mathfrak n}):=\Gamma(1)/\Gamma({\mathfrak n}) {\cal Z}$ induce
automorphisms of $X({\mathfrak n})$. It is even known that $G({\mathfrak n})$ is the full 
automorphism group of $X({\mathfrak n})$ if $p \neq 2, q \neq 3$ (cf.\ \cite{Cornelissen:01}, prop.\ 4).
It follows from (\ref{char0}) that a ``classical'' modular curve does
not admit equivariant deformations, since
$X(n) \rightarrow X(1)={\bf P}^1$ is ramified above three points (the most famous such curve is probably $X(7)$, 
which is isomorphic to Klein's quartic of genus $3$ with $PSL(2,7)$
as automorphism group). What is the analogous
result for the Drinfeld modular curves $X({\mathfrak n})$? Note that 
$X({\mathfrak n}) \rightarrow X(1)={\bf P}^1$ is ramified above two points
with ramification groups $\Z/(q+1)\Z$ and $(\Z/p\Z)^{td} \sd \Z/(q-1)\Z$
where $d=\deg({\mathfrak n})$ (cf.\ loc.\ cit.), which already shows that something
goes wrong when dully applying (\ref{char0}).

\medskip

The correct framework for carrying out such computations is that of
equivariant cohomology (Grothendieck \cite{Gr57}). It predicts that
in positive characteristic, the group cohomology of the ramification 
groups with values in the tangent sheaf at branch points contribute
to the deformation space. Bertin and M{\'e}zard have considered the case
of a cyclic group of prime order $G=\Z/p\Z$ in \cite{BM00}, mainly
concentrating on mixed characteristic lifting. In this paper, we work in 
the equicharacteristic case and
we do not want to impose direct restrictions on the group $G$, but
rather on the curve $X$, which we require to be {\sl ordinary.} 
This means that the $p$-rank of its Jacobian satisfies
$$ \dim_{{\bf F}_p} \mbox{Jac}(X)[p] = g. $$
The property of being ordinary is open and dense in the moduli space 
of curves of genus $g$, so our calculations do apply to a ``large portion'' 
of that moduli space (there is some cheating here, since curves 
without automorphisms are also dense; but notice that the analytic 
constructions in part B of the paper at least show the existence
of lots of ordinary curves with automorphisms). 
The main advantage of working with ordinary curves is that their
ramification groups are of a very specific form, so the Galois cohomological 
computation becomes feasible (cf.\ prop.\ \ref{para-ordinary-ramification}): if $P_i$ is a point on $Y$ branched
in $X \rightarrow Y$, then 

\begin{equation} \label{rampoints} 
G_{P_i} \cong (\Z/p\Z)^{t_i} \sd \Z/n_i\Z
\end{equation} 

\noindent for
some integers $(t_i,n_i)$ satisfying $n_i|p^{t_i}-1$. If $t_i>0$, we
say that $P_i$ is wildly branched.
Our main algebraic result is the following (which we state here
in a form which excludes a few anomalous cases); one can think of it giving the ``positive characteristic error term'' to the Riemann surface
computation in (\ref{char0}).

\medskip

{\bf Main Algebraic Theorem} (cf.\ \ref{global-theorem}). \ {\sl Assume that 
$p \neq 2,3$, $X$ is an ordinary curve of genus $g \geq 2$ over a field of characteristic $p>0$, 
$G$ is a finite group acting via $\rho : G \rightarrow \mbox{Aut}(X)$ on 
$X$, such that $X \rightarrow Y:=G \backslash X$ is branched
above $n$ points, of which the first $s$ are wildly branched. Then 
the equicharacteristic deformation functor $D_{X,\rho}$ is pro-representable by a ring
$H_{X,\rho}$ whose Krull-dimension is given by 
$$ \dim H_{X,\rho} = 3 g_Y -3 + n + \sum_{i=1}^s \frac{t_i}{s(n_i)}, $$
where $s(n_i):=[\F_p(\zeta_{n_i}):\F_p]=\min \{ s'>0:n_i|p^{s'}-1 \}$, $g_Y$ is the genus of $Y$ 
and $(t_i,n_i)$ are the data corresponding to the wild 
ramification points via \textup{(\ref{rampoints})}.}

\medskip

The proof goes by first computing the first-order local deformation 
functors, then studying their liftings to all of $\C_k$, and putting
the results together via the localization theorem of \cite{BM00}. In the 
end, we can even describe the ring $H_{X,\rho}$ explicitly (cf. \ref{pro-globalhull} and \ref{hull}). It turns out
to be a formal polydisc to which for each ramification point with $n_i \leq 2$
a $\frac{p-1}{2}$-nilpotent 0-dimensional
scheme is attached (and in characteristic $p=2$ even more weird things can happen). These nilpotent schemes are related to the lifting of a specific
first order deformation via what we call {\sl formal truncated Chebyshev 
polynomials}.

\medskip

{\bf Example} (continued). \ For the Drinfeld modular curve $X({\mathfrak n})$, 
we find a $(d-1)$-dimensional reduced deformation space. 

\medskip

The second part of this paper is concerned with analytic equivariant deformation theory. Let $(K,|\cdot|)$ be a non-archimedean valued field of 
positive characteristic $p>0$ with residue field $k$. Recall that a {\sl Mumford curve} over $K$
is a curve $X$ whose stable reduction is isomorphic to a union 
of rational curves intersecting in $k$-rational points. Mumford has 
shown that this is equivalent to its analytification $X^\an$ being isomorphic to an analytic space of the form $\Gamma \backslash ({\bf P}^{1,\an}_K-{\cal L}_\Gamma)$, where $\Gamma$ is a discontinuous group in $PGL(2,K)$ with ${\cal L}_\Gamma$ as set of limit points. It is known that Mumford curves are ordinary, 
so the above algebraic theory applies. But what interests us most
in this second part is to find out where these deformations ``live''
in the realm of discrete groups. 

The set-up is as follows: let $X$ be a Mumford curve, and let $N$ be a
group contained in the normalizer $N(\Gamma)$ of the corresponding
``Schottky'' group $\Gamma$ such that $\Gamma \subseteq N$. Then there
is an injection $\rho: G:= N/\Gamma \hookrightarrow \A(X)$ (isomorphism if $N=N(\Gamma)$). By rigidity, two Mumford curves are isomorphic if
and only if their Schottky groups are conjugate in $PGL(2,K)$. Hence it
is natural to consider the analytic deformation 
functor
$D_{N,\phi} : \C_K \rightarrow  {\tt Set}$ that associates to $A \in \C_K$ the set 
of homomorphisms $N \rightarrow PGL(2,A)$ that lift the given morphism
$\phi$ in the obvious sense. This functor comes with a natural
action of conjugation by the group functor $PGL(2)^\wedge : \C_K \rightarrow \mbox{\tt{Groups}}$ given by 
$$ PGL(2)^\wedge(A) := \ker [PGL(2,A) \rightarrow PGL(2,K)], $$ and we denote by $D^\sim_{N,\phi}$ the quotient functor 
$D^\sim_{N,\phi} = PGL(2)^\wedge \backslash D_{N,\phi}$.  We
now set out to compute the hulls corresponding to these functors. 

For this, we use the fact that $N$ can be described by the theorem of
Bass-Serre, which states that there is a graph of groups $(T_N,N_\bullet)$
such that $N$ is a semidirect product of a free group (of rank the
cyclomatic 
number of $T_N$) and the tree (``amalgamation'') product associated to 
a lifting of stabilizer groups $N_\ast$ of vertices and edges $\ast$ from $T_N$ to the universal covering of
$T_N$ as as graph. In a technical proposition, we show that this
decomposition of the group $N$
induces a decomposition of functors
$$ D_{N,\phi} \cong \lim_{\stackrel{\longleftarrow}{T_N}}
D_{N_\bullet,\phi|_{N_\bullet}} \times D_{F_c,\phi \circ s}, $$
where $s$ is a section of the semi-direct product (note that this is
a {\sl direct product} of functors, the inverse limit is over
$T_N$, {\sl not} over
its universal covering, and that we are using $D_{N,\phi}$, not $D^\sim_{N,\phi}$).
This reduces everything to the computation of $D_{N,\phi}$ for free
$N$ (which is easy) or finite $N$ acting on ${\bf P}^1$. The latter can be done
by using the classification of finite subgroups of ${\bf P}^1$ and the
algebraic
results from the first part in the particular case of ${\bf P}^1$. Modulo a few
anomalous 
cases, the result is the following:

\medskip

{\bf Main Analytic Theorem} (cf.\ \ref{main-analytic}). \ {\sl If $X$ is a Mumford curve of 
genus $g \geq 2$ over a non-archimedean field $K$ of characteristic
$p>3$ with Schottky group $\Gamma$, then for a given discrete
group $N$ contained in the normalizer of $\Gamma$ in $PGL(2,K)$ with corresponding
graph of groups $(T_N,N_\bullet)$, the equicharacteristic analytic deformation  
functor $D^\sim_{N,\phi}$ is pro-representable by a ring $H^\sim_{N,\phi}$ whose
dimension satisfies 
$$ \dim H^\sim_{N,\phi} = 3 c(T_N) - 3 + \sum_{{v = \mbox{{\tiny vertex}}\atop{\mbox{{\tiny of }} T_N}}} d(v) - \sum_{{e = \mbox{{\tiny edge}}\atop{\mbox{{\tiny of }} T_N}}} d(e), $$
where $c(\ast)$ denotes the cyclomatic number of a graph $\ast$ and
$$ d(\ast) = \left\{ \begin{array}{lll} 
2        & \mbox{ if } & N_\ast = \Z/n\Z , \\
3        & \mbox{ if } & N_\ast = A_4, S_4, A_5, D_n, PGL(2,p^t), PSL(2,p^t) , \\
t        & \mbox{ if } & N_\ast = (\Z/p\Z)^t,\\
t/s(n)+2 & \mbox{ if } & N_\ast =(\Z/p\Z)^t \sd \Z/n\Z, \end{array} \right. $$
where $n$ is coprime to $p$.}

\medskip

{\bf Example} (continued). \ The Drinfeld modular curve $X({\mathfrak n})$ is known
to be a Mumford curve, and the normalizer of its Schottky group is
isomorphic to an amalgam (cf.\ \cite{Cornelissen:01}) 

\begin{equation} N({\mathfrak n}) =  PGL(2,q) \ast_{({\bf Z}/p\Z)^t \sd
    {\bf Z}/{(q-1)\Z}} ({\bf Z}/p\Z)^{td} \sd {\bf Z}/{(q-1)}\Z.
  \label{NN1} \end{equation}

\noindent The above formula gives again a $(d-1)$-dimensional deformation space. We can see these deformations explicitly as follows: by conjugation
with $PGL(2,C)$, we can assume that the embedding of $PGL(2,p^t)$ in (\ref{NN1}) into $PGL(2,C)$ is 
induced by the standard $\F_{q} \subseteq C$. Then a little matrix
computation shows that the freedom of choice left is that of a $d$-dimensional
$\F_{q}$-vector space $V$ of dimension $d$ in $C$ which contains the standard $\F_q$ in order to embed the order-$p$ elements from the second group involved in the amalgam as
${1 \, x}\choose{0 \, 1}$ for $x \in V$. 

\medskip

In the final section of the paper, we compare algebraic and analytic
deformation functors. What comes out is an isomorphism of functors
$$ D^\sim_{N,\phi} \cong D_{X,\rho}.$$ To prove this, we remark that the whole construction
of Mumford can be carried out in the category $\C_K$ (by fixing a lifting of $N$ to $PGL(2,A)$) to produce a scheme
$X_\Gamma$ over $\Spec \ A$ whose central fiber is isomorphic
to $X$, and such that $X_\Gamma$ carries an action $\rho$ of $N/\Gamma$ which
reduces to the given action on $X$. 

\medskip

{\bf Remark.} \ We did not take the more ``global'' road of
``equivariant Teichm{\"u}ller space'' to
analytic deformation, as Herrlich does in \cite{Herrlich:84} and
\cite{Herrlich:87}: one can consider the space 
$$ \M(N,\Gamma) = \PGL(2,K) \backslash \Hom^*(N,PGL(2,K)) / \Aut_\Gamma(N),$$
where $\Hom^*$ means the space of injective morphisms with discrete
image, the action of $PGL(2,K)$ from the left is by conjugation, and
the 
action of the group of automorphisms of $N$ which fix $\Gamma$ is on
the right. The relations in $N$ impose a natural structure on
$\M(N,\Gamma)$ as an analytic space, but in positive characteristic,
this structure might be non-reduced due to the presence of parabolic
elements, 
e.g., if $N$ contains $\Z/p\Z, p>3$, then the condition that a matrix
$\gamma = {{a\,b}\choose{c\,d}}$ be of order $p$ leads to 
$$ (\mathrm{tr}^2(\gamma)-4 \det(\gamma))^{\frac{p-1}{2}} \cdot b = 0. $$ 
Therefore in \cite{Herrlich:87}, global considerations are
restricted
to the case where $N$ does not contain parabolic elements (loc.\ cit., p.\ 148),
whereas our local calculations are independent of such restrictions. 

Let us note that if $K$ has characteristic 0, then  in \cite{Herrlich:87}
a formula is given for the dimension of $\M(N,\Gamma)$ (which turns out to be an equidimensional space) compatible with our main analytical theorem.

\medskip

{\bf Convention.} \ Throughout this paper, if $R$ is a local ring, we
let $\dim R$ denote its Krull-dimension, and if $V$ is a $k$-vector space,
we let $\dim_k V$ denote its dimension as a $k$-vector space.

\vspace{3ex}

{\bf PART A: ALGEBRAIC THEORY}

\vspace{2ex}

\sectioning\label{section-covering}\ 
{\bf Deformation of Galois covers}

\vspace{1ex}
\paragraph\label{global-def} {\bf The global deformation functor.} \ We will start by recalling the definition
of equivariant deformation of a curve with a given group of automorphisms
(cf.\ \cite{BM00} for an excellent
survey).

Let $k$ be a field, and $X$ a smooth projective curve over $k$. We fix a 
finite subgroup $G$ in the automorphism group $\Aut_k(X)$ of $X$, and denote
by $\rho$ the inclusion $$\rho  :  G\hookrightarrow\Aut_k(X) : \sigma \rightarrow \rho_\sigma.$$
Let $Y=G\backslash X$ be the quotient curve, and denote by $\pi$ the quotient map $X\rightarrow Y$.

Let $\C_k$ be the category of Artinian local $k$-algebras with the residue 
field $k$. 
A {\sl lifting} of $(X,\rho)$
to $A$ is a triple $(X^{\sim},\rho^{\sim},\phi^{\sim})$ consisting of a scheme $X^{\sim}$ which is smooth of finite 
type over $A$, 
an injective group homomorphism $$\rho^{\sim}\colon G\hookrightarrow
\Aut_A(X^{\sim}) : \sigma \rightarrow \rho_\sigma^\sim,$$ and an isomorphism
$$\phi^{\sim}\colon X^{\sim}\otimes_{\, \Spec \, A}\Spec  \ k\rightarrow X$$ of 
schemes over 
$k$ such that $\ovl{\rho^{\sim}}=\rho.$ Here, $\ovl{\rho^{\sim}}$ denotes the
composite of $\rho^{\sim}$ and restriction onto the central
fiber, identified with $X$ by $\phi^{\sim}$.

Two liftings $(X^{\sim},\rho^{\sim},\phi^{\sim})$ and 
$(X^{\approx},\rho^{\approx},\phi^{\approx})$ are said to be 
{\sl isomorphic} if there exists an isomorphism
$\psi\colon X^{\approx}\rightarrow X^{\sim}$ of schemes over $A$ such that
$\phi^{\sim}\circ(\psi\otimes_Ak)=\phi^{\approx}$ and, for any $\sigma\in G$, $\psi\circ\rho^{\approx}_{\sigma}=
\rho^{\sim}_{\sigma}\circ\psi$.

We arrive at the {\sl deformation functor} 
$$ D_{X,\rho} : \C_k \rightarrow {\tt Set}$$ that assigns to 
any $A \in \C_k$ the set of isomorphism classes of liftings of $(X,\rho)$. 

\vspace{1ex}
\paragraph\label{para-tangentmap} {\bf The functor  $\pi^G_{\ast}$ on tangents.} \ 
In order to compute the tangent space to this deformation functor, we
need to recall the equivariant cohomology theory of Grothendieck.

The morphism $\pi$ induces a tangent map $\T_X\rightarrow\pi^{\ast}\T_Y$, 
which is a monomorphism since $\pi$ is generically {\'e}tale (dually, $\Omega_{X/Y}$ is a torsion sheaf).
Note that both $\T_X$ and $\pi^\ast \T_Y$ are $G$-$\O_X$-modules
(cf.\ \cite[5.1]{Gr57}), and the tangent map $\T_X\rightarrow\pi^{\ast}\T_Y$
is a morphism of $G$-$\O_X$-modules; the $G$-structure of 
$\pi^\ast \T_Y=\O_X\bigotimes_{\O_Y}\pi^{-1}\T_Y$ is the tensor product of the usual 
$G$-structure on $\O_X$ and the trivial $G$-structure on $\pi^{-1}\T_Y$
(cf.\ [loc.\ cit.]).
We apply the functor $\pi^G_{\ast}$ to the tangent map; recall that by definition, for a $G$-sheaf $\cal F$ on $X$ and an open $U$ of $Y$, $\pi^G_\ast({\cal F})(U)$ is the set of sections of $\cal F$ on $\pi^{-1}(U)$
invariant under $G$. Since 
$\pi^G_{\ast}\pi^{\ast}\T_Y\cong\T_Y$ (cf.\ \cite[(5.1.1)]{Gr57}), we get
$$
0 \rightarrow \pi^G_{\ast}\T_X\longrightarrow\T_Y,
\leqno{\indent\textrm{{\rm (\ref{para-tangentmap}.1)}}}
$$
The fact that this map is a monomorphism is due to the left exactness of $\pi^G_{\ast}$. 
Note that $\pi^G_{\ast}\T_X$ is an invertible sheaf on $Y$: indeed, it can be
described by formula (\ref{pro-linebundle}.1) below.

\vspace{1ex}
\paragraph\label{para-ordinary} {\bf Ordinary curves.} \ 
From now on, we will assume that the field $k$ is of characteristic 
$p>0$ and that the curve $X$ is {\sl ordinary}. This means that
$$ \dim_{{\bf F}_p} \mbox{Jac}(X)[p] = g. $$
The property of being ordinary is open and dense on the moduli space 
of curves of genus $g$ (\cite{Nakajima:87}, \S 4, \cite{Koblitz:75}), since
it is essentially a maximal rank condition on the Hasse-Witt matrix of the 
curve (hence open), and ordinary curves exist for all $g$. 

\vspace{1ex}
\paragraph\label{para-ordinary-ramification}\ 
{\bf Proposition.} {\sl Let $X$ be an ordinary curve and $G$ a finite group
of automorphisms of $X$, $\pi : X \rightarrow G \backslash X$ and $P \in Y$ a branch point of $\pi$. Then the ramification filtration for $P$ stops at
$G_2 = \{ 1 \}$. The ramification group at $P$ is of the form $({\Z}/p\Z)^t \sd \Z/n\Z$
for $n|p^t-1$. Actually, letting $s:=[\F_p(\zeta):\F_p]$ where $\zeta$ is a primitive $n$-th root of unity, one can consider $({\Z}/p\Z)^t$
as a vector space of dimension $\frac{t}{s}$  over $\F_q$ where $q=p^s$, and the action 
of $\Z/n \Z$ is exactly by multiplication with $\zeta$. Locally at $P$, the
cover can be decomposed as a tower $$k((x_{\frac{t}{s}})) - \ldots - k((x_2))- k((x_1)) - k((x_0))$$ where the first step is a Kummer extension of degree $n$ ($x_0=x_1^n$) and all others are Artin-Schreier extensions of degree $q$
($x_i^{-q}-x_i^{-1}=x_{i-1}^{-1}$).} 

\vspace{1ex}
\pf
Let $P$ be a ramified point in that cover and let $G_0 \supseteq G_1 \supseteq G_2 \supseteq \dots$ be the standard filtration on the ramification group
$G_0$ at $P$. In \cite{Nakajima:87} (Theorem 2(i)), S.\ Nakajima has shown  that           for an ordinary curve, $G_2= \{ 1 \}$. It is also known that $G_i/G_{i+1}$ for $i \geq 1$ are elementary abelian $p$-groups and that $G_0$ is the semi-direct 
product of a group of order prime to $p$ and $G_1$ with the action given above (\cite{Serre:68}, IV.2 cor.\ 1-4). The results follows from this and standard field theory.
\qed

\vspace{1ex}
{\bf \ref{para-ordinary-ramification}.1 Remark.} \  The following formula for $s$  holds ((\cite{Lidl}, 2.47(ii)): 
$$s = \min \{ s'  \in \Z_{>0} :   n|p^{s'}-1 \}.$$ 
 
\vspace{1ex}
{\bf \ref{para-ordinary-ramification}.2 Remark.} \ 
 The calculations that follow actually only depend upon
the special form of the ramification filtration given in \ref{para-ordinary-ramification}, so they also apply to (not necessarily 
ordinary) curves for which the action of $G$ is as described here.

\vspace{1ex}
\paragraph\label{para-branchpoints}\ 
{\bf Notations.}
We break up the branch points $P$ of $\pi$, which we enumerate as $\{P_1, \dots, P_n \}$ into two sets $T$ and $W$; let ${\bf Z}_p^{t} \sd {\bf Z}_{n}$ be the ramification group at $P$. 
Then 
$$ P \in T \iff t=0 \mbox{ or } (p=2 \mbox{ and } t=1), $$
and $P \in W$ otherwise. Note that for $p \neq 2$, points in $T$ are called {\sl tame} and points in $W$ are called {\sl wild}, but for $p=2$ the distinction is
more subtle. Observe that $p=2$ and $t=1$ implies $n=1$.
We define the divisor $\Delta$ on $Y$ by
$$
\Delta=
{\textstyle \sum\limits_{P \in T} P +2\sum\limits_{Q \in W} Q} \ \ \  \mbox{ of degree }  \ \ \ \delta:=\deg(\Delta) = |T|+2|W|.$$

\vspace{1ex}

\paragraph\label{pro-linebundle}\
{\bf Proposition.}\  $\pi^G_{\ast}\T_X\cong\T_Y(-\Delta)$.

\vspace{1ex}
\pf
Like in \cite{BM00} (proof of 5.3.2), we have that 
$$\pi^G_{\ast}\T_X\cong \T_Y \otimes [ \O_Y \cap \pi_\ast(\O_X(-{\mathfrak r}))], 
\leqno{\indent\textrm{{\rm (\ref{pro-linebundle}.1)}}}
$$
where $\mathfrak r$ is the (global) different of $k(X)/k(Y)$. The divisor $\mathfrak{r}$ is supported at branch points $P$. The exponent of the local different at $P$ is $\sum_{i=0}^\infty (|G_i|-1)$ where $G_i$ are the higher ramification groups at $P$. Using \ref{para-ordinary-ramification} one gets that $ {\mathfrak r}_P =  (np^t-1+p^t-1) \cdot P$.
Hence upon intersecting with $ \O_Y$ in (\ref{pro-linebundle}.1), 
we get that (recall that the cover is totally ramified)
\begin{eqnarray*}
[ \O_Y \cap \pi_\ast(\O_X(-{\mathfrak r}))]_P  &=& \O(- \lceil 1 +  \frac{p^t-2}{np^t} \rceil P) \\
&=& \left \{ \begin{array}{ll} \O(-P) & \mbox{if } t=0 \mbox{ or } (p=2,t=1); \\ \O(-2 P) & \mbox{otherwise.} \end{array} \right. \end{eqnarray*}
where $\lceil x \rceil = \min \{ n \in {\bf Z}_{>0} \ : \  n \geq x \}$. If we now
collect the local terms, the result comes out. \qed

\vspace{1ex}
{\bf \ref{pro-linebundle}.2 Remark.} \ On pages 235-236 of \cite{BM00}, $\lceil \cdot \rceil$ and 
$\lfloor \cdot \rfloor$ have to be interchanged everywhere.

\vspace{1ex}
\paragraph\label{para-spectral} {\bf Cohomology of $\pi_\ast^G$.} \ 
We now recall the two spectral sequences from \cite[5.2]{Gr57}.
Let ${\cal F}$ be a coherent $G$-$\O_X$-module, and write
$$
\begin{array}{lcl}
\H^i(X;G,{\cal F})&=&\mathrm{R}^i\Gamma^G_X{\cal F}, \\
{\cal H}^i(G,{\cal F})&=&\mathrm{R}^i\pi^G_{\ast}{\cal F},
\end{array}
$$
where $\Gamma^G_X=\Gamma_Y\circ\pi^G_{\ast}$, i.e., $\Gamma^G_X{\cal F}=
(\Gamma_X{\cal F})^G$.
There are two cohomological spectral sequences
$$
\begin{array}{lcl}
\llap{${}^I\!E^{p,q}_2=$\ }\H^p(Y,{\cal H}^q(G,{\cal F}))&\Longrightarrow&
\H^{p+q}(X;G,{\cal F}), \\
\llap{${}^{II}\!E^{p,q}_2=$\ }\H^p(G,\H^q(X,{\cal F}))&\Longrightarrow&
\H^{p+q}(X;G,{\cal F}).
\end{array}
$$
The first one gives rise to the edge sequence
$$
0\longrightarrow\H^1(Y,\pi^G_{\ast}{\cal F})\longrightarrow
\H^1(X;G,{\cal F})\longrightarrow\H^0(Y,{\cal H}^1(G,{\cal F}))
\longrightarrow 0
\leqno{\indent\textrm{{\rm (\ref{para-spectral}.1)}}}
$$
as we are on a curve $Y$.

\vspace{1ex}

\paragraph\label{globaltangentintro} {\bf Tangent space to the global deformation functor.} \ 
We can use this equivariant cohomology to compute the tangent space
to our deformation functor. Recall that the ring of dual numbers
is defined as $k[\epsilon]=k[E]/(E^2)$; clearly, it belongs to $\C_k$.
The tangent space to the deformation functor $D_{X,\rho}$ is by definition
its value on the ring of dual numbers with its natural $k$-linear
structure (cf.\ \cite{Sc68}).

\vspace{1ex}

\paragraph\label{globaltangent}\
{\bf Proposition} (\cite{BM00}, (3.2.1), (3.3.1)) \ {\sl We have $D_{X,\rho}(k[\epsilon])\cong
\H^1(X;G,\T_X)$. \qed}

\vspace{1ex}
\paragraph\label{para-globallocal}\ {\bf Localization.} \ 
We will now describe localization for our deformation functor. Let $Q_j\in X$ be a point lying over $P_j$ for $1\leq j\leq n$, and $G_j$ the
stabilizer of $Q_j$, which acts on the local ring $\O_{X,Q_j}$ by 
$k$-algebra automorphisms; we denote it by $\rho_j\colon G_j\rightarrow
\Aut_k(\O_{X,Q_j})$. Changing the choice of $Q_j$ does not affect, up to
a suitable equivalence, on the action $\rho_j$.
Every lifting $(X^{\sim},\rho^{\sim},\phi^{\sim})$ induces a {\sl lifting
of the local representation $\rho_j$} for any $1\leq j\leq n$ in the 
following sense: 

\vspace{1ex}
\paragraph\label{para-globallocal2}\ {\bf The local deformation functor.} \ 
Let $R$ be a $k$-algebra, and $G$ a group acting via $\rho : G \rightarrow \Aut_k(R)$ on $R$ by $k$-algebra
automorphisms.  Let $A \in \C_k$ and  $R_A = R \otimes_k A$. Let  $\ovl{\,\cdot\,}$ denote the reduction 
map $R_A \rightarrow R$ modulo the maximal ideal $m_A$ of $A$.

We define a lifting of $\rho$ to $A$ as
a group homomorphism $\rho^{\sim}\colon G\rightarrow\Aut_{A}(R_A)$ 
such that $\ovl{\rho^{\sim}}=\rho$. 
Two liftings $\rho^{\sim}$ and $\rho^{\approx}$ are said to be 
{\sl isomorphic} (over $\rho$) if there exists 
$\psi\in\Aut_A(R_A)$ with $\ovl{\psi}=\mathrm{Id}_R$
such that, for any $\sigma\in G$, $\psi\circ\rho^{\approx}_{\sigma}=
\rho^{\sim}_{\sigma}\circ\psi$.
We arrive at a {\sl local deformation functor} $$
D_{\rho}\colon\C_k\longrightarrow{\tt Set}
$$
such that $D_\rho (A)$ is the set of all isomorphism classes of liftings of 
$\rho$ to $A$. It is but a formal check using Schlessinger's criterion 
to see that $D_{\rho}$ has a pro-representable hull (which we will denote
by $H_{\rho}$) if $G$ is a finite group (compare \cite{BM00}, 2.2). 

\vspace{1ex}
 
In our situation, we get a transformation of functors
$$
D_{X,\rho}\longrightarrow D_{\rho_1}\times\cdots\times D_{\rho_n}
\leqno{\indent\textrm{{\rm (\ref{para-globallocal}.1)}}}
$$
It turns out that this morphism is formally smooth (\cite{BM00}, 3.3.4). Hence
we get the following localization result:

\vspace{1ex}

\paragraph\label{pro-globalhull}\
{\bf Proposition} (\cite{BM00}, 3.3.1, 3.3.5).\ {\sl 
The functor $D_{X,\rho}$ has a pro-representable hull $H_{X,\rho}$; in fact
$$
H_{X,\rho}\cong H_{\rho_1}\hat{\otimes}\cdots\hat{\otimes}H_{\rho_n}
[[u_1,\ldots,u_N]],
$$
where $N=\dim_k\H^1(Y,\pi^G_{\ast}\T_X)$. \qed}

\vspace{1ex}

\paragraph\label{representability}\
{\bf Remark.} \ As soon as $g \geq 2$, the curve $X$ does not
have infinitesimal automorphisms, so the functor $D_{X,\rho}$ is in 
fact {\sl pro-representable} by $H_{X,\rho}$ (cf.\ \cite{BM00}, 2.1).

\vspace{2ex}

\sectioning\label{section-lifting}\ 
{\bf Lifting of group actions and group cohomology}

\vspace{1ex}

The aim of the next three sections is to compute the tangent space and the pro-representable hull of the local deformation functors for 
the representations of the branch groups in quotients of 
ordinary curves as automorphisms of the local stalks of the tangent sheaf. 

\vspace{1ex}
\paragraph\label{para-tangent} {\bf Action on derivations.} \ 
Let $k$ be a field and $R$ a $k$-algebra. We denote by $\T_R$ the $R$-module 
of $k$-derivations, i.e., the set of all maps $\delta\colon R\rightarrow R$
such that $\delta(xy)=\delta(x)y+x\delta(y)$ for $x,y\in R$ and that
$\delta(a)=0$ for $a\in k$.
Let $\varphi\in\Aut_k(R)$ be an automorphism of $R$ over $k$.
Then it induces a $k$-module automorphism of $\T_R$, denoted by 
$\Ad_{\varphi}$, by
$$
\Ad_{\varphi}(\delta)=\varphi\circ\delta\circ\varphi^{-1}.
$$
Note that $\Ad_{\varphi}$ is not an $R$-module automorphism, but rather equivariant:  for $x\in R$ and
$\delta\in\T_R$, we have $\Ad_{\varphi}(x\delta)=\varphi(x)\Ad_{\varphi}(\delta)$. Let $G$ be a finite group and suppose it acts on $R$ by $k$-algebra automorphisms; 
i.e., suppose a group homomorphism $$\rho\colon G\rightarrow\Aut_k(R) \ : \ 
\sigma\mapsto\rho_{\sigma} $$ is given.
The induced action of $G$ on $\T_R$ by $k$-module automorphisms is denoted
by $$\Ad_{\rho} \ : \ \sigma\mapsto\Ad_{\rho,\sigma}. $$

\vspace{1ex}
\paragraph\label{para-dualnumber} {\bf Tangent space to the local deformation functor.} \ 
Let $k[\epsilon]$ be the ring of dual numbers.
Then $R[\epsilon]=R\bigotimes_kk[\epsilon]$ is the 
$k[\epsilon]$-algebra of elements $x+y\varepsilon$ with
$x,y\in R$. Let $\ovl{\,\cdot\,} \ : R[\epsilon] \rightarrow R$ denote the  reduction map modulo-$\epsilon$.
In the following proposition, we make an explicit identification
between the tangent space  $D_\rho(k[\epsilon])$ to the deformation functor
$D_\rho$ and the first group cohomology with values in the derivations. 

\vspace{1ex}

\paragraph\label{pro-lifting}\
{\bf Proposition.}\  {\sl 
There exists a bijection \textup{(}depending on the deformation parameter
$\epsilon$\textup{)}
$
d \ \colon \ D_{\rho}(k[\epsilon])
\stackrel{\sim}{\longrightarrow}
\H^1(G,\T_R)
$
described as follows: The $1$-cocycle $d_{\rho^{\sim}}$ associated to a
lifting $\rho^{\sim}$ is given by the formula
$$
d_{\rho^{\sim}} \sigma=\frac{\rho^{\sim}_{\sigma}\circ\rho^{-1}_{\sigma}
-\mathrm{Id}}{\epsilon}\ 
 \ (={\frac{d}{d\epsilon}}(\rho^{\sim}_{\sigma}\circ\rho^{-1}_{\sigma})
|_{\epsilon=0}), 
$$
for any $\sigma\in G$.}

\vspace{1ex}
\pf For $\sigma\in G$ we set $\rho^{\sim}_{\sigma}(x)=\rho_{\sigma}(x)
+\rho'_{\sigma}(x)\epsilon$ ($x\in R$).
Then for $x+y\epsilon\in R[\epsilon]$, we have
$$
\rho^{\sim}_{\sigma}(x+y\epsilon)=\rho^{\sim}_{\sigma}(x)+
\rho^{\sim}_{\sigma}(y)\epsilon=
\rho_{\sigma}(x)+[\rho'_{\sigma}(x)+\rho_{\sigma}(y)]\epsilon,
$$
i.e., $\rho'_{\sigma}$ determines the lifting $\rho^{\sim}$.
The $1$-cocycle $d_{\rho^{\sim}}$ is given by $d_{\rho^{\sim}}\sigma=
\rho'_{\sigma}\circ\rho^{-1}_{\sigma}$.
The following two formulas are straightforward:
\begin{condlist}
\item[(i)] $\rho'_{\sigma}(xy)=\rho_{\sigma}(x)\rho'_{\sigma}(y)+
\rho'_{\sigma}(x)\rho_{\sigma}(y)$ for $x,y\in R$.
\item[(ii)] $\rho'_{\sigma\tau}=\rho'_{\sigma}\circ\rho_{\tau}+
\rho_{\sigma}\circ\rho'_{\tau}$ for $\sigma,\tau\in G$.
\end{condlist}
From (i) it follows that $d_{\rho^{\sim}}\sigma\in\T_R$, and from (ii) 
that $d_{\rho^{\sim}}$ is a cocycle, i.e., 
$d_{\rho^{\sim}}\sigma\tau=d_{\rho^{\sim}}\sigma+\Ad_{\rho,\sigma}
(d_{\rho^{\sim}}\tau)$.

Suppose two liftings $\rho^{\sim}$ and $\rho^{\approx}$ are isomorphic by
$\psi\in\Aut_{k[\epsilon]}(R[\epsilon])$. By the condition 
$\ovl{\psi}=\mathrm{Id}_R$ we can write $\psi(x)=x+\delta(x)\epsilon$ 
($x\in R$), where $\delta\in\T_R$.
The equality $\psi\circ\rho^{\approx}_{\sigma}=
\rho^{\sim}_{\sigma}\circ\psi$ implies that 
$$\rho''_{\sigma}-\rho'_{\sigma}=
\rho_{\sigma}\circ\delta-\delta\circ\rho_{\sigma},
\leqno{\indent\textrm{{\rm (\ref{pro-lifting}.1)}}}
$$ where
$\rho''_{\sigma}$ is the $\epsilon$-part in $\rho^\approx$ (as $\rho'_{\sigma}$ was the $\epsilon$-part of $\rho^{\sim}$).
Hence $$d_{\rho^{\approx}}\sigma-d_{\rho^{\sim}}\sigma=
\Ad_{\rho,\sigma}(\delta)-\delta,
\leqno{\indent\textrm{{\rm (\ref{pro-lifting}.2)}}}
$$ which implies that $(d_{\rho^{\approx}})-
(d_{\rho^{\sim}})$ is a coboundary.
Conversely, if we have (\ref{pro-lifting}.2) for some $\delta\in\T_R$, then one can
define $\psi\in\Aut_{k[\epsilon]}(R[\epsilon])$ by the obvious formula
$\psi(x)=x+\delta(x)\epsilon$ (or, equivalently, $\psi(x+y\epsilon)=
x+[y+\delta(x)]\epsilon$), which gives an isomorphism between the liftings
$\rho^{\sim}$ and $\rho^{\approx}$.
Therefore, the map $d$ is well-defined, and is injective.

We now show surjectivity.
A given cocycle $d \ \colon \ G\rightarrow\T_R$ induces an automorphism $\rho^{\sim}_{\sigma}$ for any $\sigma\in G$
by the formula: $\rho^{\sim}_{\sigma}(x)=\rho_{\sigma}(x)+(d\sigma)\circ
\rho_{\sigma}(x)\epsilon$ for $x\in R$.
One can easily check that $\sigma\mapsto\rho^{\sim}_{\sigma}$ gives a
lifting of $\rho$ whose associated $1$-cocycle is exactly
$d$.
\qed

\vspace{2ex}
\sectioning\label{section-groupcohomology}\ 
{\bf Computation of group cohomology}

\vspace{1ex}
\paragraph\label{para-cases}\ 
In this section, $k$ denotes a field of characteristic $p>0$, and $\O$ a 
discrete 
valuation ring over $k$ with the residue field $k$. We fix a regular parameter
$x$ for $\O$. 
The $\O$-module $\T_{\O}$ (defined in (\ref{para-tangent})) is free of rank $1$ since every $\delta\in
\T_{\O}$ is determined by $\delta(x)$. 
Let $\frac{d}{dx}$ be the unique $k$-derivation such that $\frac{d}{dx}(x)=1$.
Then $\T_{\O}=\O\frac{d}{dx}$.
Let a group $G$ act on $\O$ by $k$-algebra automorphisms. In this section 
we write the action of $\rho$ exponentially ($f\rightarrow f^{\sigma}$ for $f\in\O$ and
$\sigma\in G$), omitting $\rho$ from the notation --
we do the same for the induced action $\Ad$ on the
tangent space $\T_{\O}$.
The $G$-equivariancy 
condition now becomes $(f\delta)^{\sigma}=f^{\sigma}\delta^{\sigma}$ for
$f\in\O$.
In this section, we will compute the group cohomology $\H^1(G,\T_{\O})$ 
in the following three situations, which are exactly the ones 
that arise for the local action of the ramification group
of a branch point on an ordinary curve:

\vspace{1ex}
\indent{\rm (\ref{para-cases}.1)}\ $G=\langle\tau\rangle\cong\Z/n\Z$
with $(n,p)=1$ acting on $\O$ by $x^{\tau}=\zeta x$, where $\zeta$ is a 
primitive $n$-th root of unity.

\vspace{1ex}
\indent{\rm (\ref{para-cases}.2)}\ $G=\prod_{i=1}^t\langle\sigma_i\rangle
\cong(\Z/p\Z)^t$ acting on $\O$ by $x^{\sigma_i}=x/(1-u_ix)$ where $u_1,\ldots,
u_t\in k$ are linearly independent over $\F_p$. Let $V$ be the $t$-dimensional $\F_p$-vector subspace in $k$ spanned by
$u_1,\ldots,u_t$.
The group $G$ and its action are isomorphic to the vector group $V$, acting on $\O$ by $x^u=x/(1-ux)$ for $u\in V$; by this, 
we can (and will) pretend that $G=V$.

\vspace{1ex}
\indent{\rm (\ref{para-cases}.3)}\ $G = N \sd H$ where $N=\prod_{i=1}^{t} 
\langle\sigma_i\rangle \cong(\Z/p\Z)^t$ and $H = \langle\tau\rangle \cong \Z/n\Z$ with $n>1$ and
$n|p^t-1$. Let $\zeta$ be a primitive $n$-th root of unity in $k$ and $s=[\F_p(\zeta):\F_p]$ and $q=p^s$ as in \ref{para-ordinary-ramification}. 
Similarly to the previous case, we consider $N$ as as vector group $V$ 
of dimension $t/s$ over $\F_q$ acting on $\O$ by $x^u=x/(1-ux)$ for $u\in V$,
and the action of $H$ is scalar multiplication by $\zeta \in \F_q^*$ on $x$. 

\vspace{1ex}
\paragraph\label{para-(..1)}\ 
{\bf The case (\ref{para-cases}.1)}.\ 
Since the $G$-module $\T_\O$ is killed by $p$, but $p$ is prime to the order of $G$, all higher
group cohomology vanishes: $\H^n(G,\T_\O)=0$ for $n>0$ (\cite{Br82}, III (10.2)).

\vspace{1ex}
\paragraph\label{para-(..2)}\ 
{\bf The case (\ref{para-cases}.2)}.\ 
Since $(\frac{d}{dx})^u=(1-ux)^2\frac{d}{dx}$, 
the $G$-module $\T_{\O}$ is isomorphic to $\O$ with the $G$-action given by 
$$
f^u(x)=f\left(\frac{x}{1-ux}
\right)(1-ux)^2
\leqno{\indent\textrm{{\rm (\ref{para-(..2)}.1)}}}
$$
for $f \in \O, u \in V$.

\vspace{1ex}
\paragraph\label{lem-decompose}\
{\bf Lemma.}\ {\sl The $G$-module $\O$ with $G$-action 
(\ref{para-(..2)}.1) is isomorphic to $M\bigoplus x\O$, where $M=k\bigoplus kx
\bigoplus kx^2$ such that:

\vspace{1ex}
(1) The $G$-module structure of $M$ is given by
$$
(a_0+a_1x+a_2x^2)^u=
a_0+(a_1-2ua_0)x+(a_2-ua_1+u^2a_0)x^2,
$$
i.e., w.r.t\ the basis $\{ 1,x,x^2 \}$, 
$$
u\longleftrightarrow
\Phi(u)=
\left(
\begin{array}{ccc}
1&0&0\\
\llap{$-$}2u&1&0\\
u^2&\llap{$-$}u&1
\end{array}
\right).
$$

(2) The $G$-module structure on the second factor $x\O$ is the original $G$-action on $\O$, i.e.,
$f^u(x)=f(\frac{x}{1-ux})$  for $f(x)\in x\O$.}

\vspace{1ex}
\pf
Note that $\O=M\oplus x^3\O$ is a $G$-stable direct decomposition.
The action of $G$ on the first factor is as stated in (1).
For $x^2f(x)\in x^3\O$ (i.e., $f(x)\in x\O$), the action (\ref{para-(..2)}.1) 
gives $(x^2f(x))^u=x^2f(x/(1-ux))$.
\qed

\vspace{1ex}
\paragraph\label{para-(..2)-1}\ 
Let $G$ act on $x\O$ as in \ref{lem-decompose} (2), and let $F=\mbox{Frac}(\O)$ be the fraction field of $\O$. Since this action comes 
from that on $\O$ by ring automorphisms, it extends to $F$, respecting
the decomposition $F=x\O\oplus k[x^{-1}]$.
We have 
$$\H^1(G,x\O)\oplus\H^1(G,k[x^{-1}])=\H^1(G,F)=0$$ by standard facts (\cite{Serre:68}, X.1), and hence $\H^1(G,x\O)=0$. Thus $\H^1(G,\T_{\O})\cong\H^1(G,M)$.

We will now compute $\H^1(G,M)$. Since $G$ is a commutative $p$-group, 
the condition on a map $d : V \rightarrow M$ to be a cocycle ($d(u+v)=du+(dv)^u
$) implies
$$
\left\{ 
\begin{array}{lcc}
d(pu)= 0 \iff 
 du+(du)^u+(du)^{2u}+\cdots+(du)^{(p-1)u}=0 &\textrm{---}&\textrm{(i)}\\
d(u+v)=d(v+u) \iff  du + dv^u = dv + du^v 
&\textrm{---}&\textrm{(ii)}
\end{array}
\right.
$$
Let us write a cocycle $d$ as $du=a_0(u)+a_1(u)x+a_2(u)x^2$ for $u\in V$.

\vspace{1ex}
(\ref{para-(..2)-1}.1)
Calculating the matrix 
$1+\Phi(u)+\Phi(2u)+\cdots+\Phi((p-1)u)$, we deduce that condition (i) is
\begin{condlist}
\item[(1-a)] empty unless $p=2$ or $3$,
\item[(1-b)] equivalent to $a_0(u)=0$ if $p=3$,
\item[(1-c)] equivalent to $ua_0(u)+a_1(u)=0$ if $p=2$.
\end{condlist}

\vspace{1ex}
(\ref{para-(..2)-1}.2)
Condition (ii) is equivalent to $2ua_0(u)=2va_0(v)$ and
$u^2a_0(v)-ua_1(v)=v^2a_0(u)-va_1(u)$.
Hence
\begin{condlist}
\item[(2-a)] if $p \neq 2$, (ii) is equivalent to $a_0(u)=ua_0$ and $a_1(u)=u(a_1-ua_0)$, 
where $a_0$ and $a_1$ are constants independent of $u$,
\item[(2-b)] if $p=2$, (ii) together with (i) is equivalent to $a_0(u)=ua_0$ and
$a_1(u)=u^2a_0$, where $a_0$ is a constant independent of $u$.
\end{condlist}

\vspace{1ex}
Thus we get
$$
\mathrm{Z}^1(G,M)=\left\{
\begin{array}{ll}
\{du=ua_0+u(a_1-ua_0)x+a_2(u)x^2\}&(p\neq 2,3)\\
\{du=ua_1x+a_2(u)x^2\}&(p=3)\\
\{du=ua_0+u^2a_0x+a_2(u)x^2\}&(p=2)
\end{array}
\right.
$$
In each case, the cocycle condition is equivalent to the fact that the function $a_2$ satisfies 
$$a_2(u+v)=a_2(u)+a_2(v)+uv[(u+v)a_0-a_1]$$ (with $a_0=0$ if $p=3$ and
$a_1=0$ if $p=2$).
In particular, whenever $a_0=a_1=0$ the function $a_2$ is
$\F_p$-linear. Hence $\mathrm{Z}^1(G,M)$ contains the $t$-dimensional 
$k$-subspace $\Hom_{\F_p}(V,k)=V^{\ast}\bigotimes k$.
The dimension over $k$ of $\mathrm{Z}^1(G,M)$ is therefore $t+2$ if 
$p\neq 2,3$ or is $t+1$ otherwise.

\vspace{1ex}
\paragraph\label{para-(..2)-2}\ 
 Let $g=b_0+b_1x+b_2x^2$.
Then a coboundary is of the form
$$
g^u-g=-2ub_0x+(-ub_1+u^2b_0)x^2.
$$
We can then compute a  $k$-basis 
for $\H^1(G,\T_{\O})$ as follows: 

\vspace{1ex}
(\ref{para-(..2)-2}.1)
If $p\neq 2,3$, a non-trivial such cohomology class $[d_0]$ is given by the cocycle
$d_0$ with
$$
d_0u=-u + (u^2+u)x-(\frac{1}{3}u^3+\frac{1}{2}u^2+\frac{1}{6}u)x^2. $$
The other cohomology classes come from the subspace 
$\{a_0=a_1=0\}\cong\Hom_{\F_p}(V,k)$ in $\mathrm{Z}^1(G,M)$,
in which the coboundary classes are spanned by $du=ux^2$. Hence the part of the cohomology coming from $\{a_0=a_1=0\}$ is isomorphic
to $$\Hom_{\F_p}(V,k)/k\cdot\iota,$$ where 
$\iota\colon V\hookrightarrow k $ is the natural inclusion.
Hence 
$$
\H^1(G,\T_{\O})\cong k\cdot[d_0] \ {\textstyle \bigoplus} \ 
\Hom_{\F_p}(V,k)/k\cdot\iota.
$$
In particular, $\dim_k\H^1(G,\T_{\O})=t$.

\vspace{1ex}
(\ref{para-(..2)-2}.2)
If $p=3$, then only the cocycles coming from $\{a_0=a_1=0\}$ survive.
Hence 
$$
\H^1(G,\T_{\O})\cong \Hom_{\F_3}(V,k)/k\cdot\iota,
$$
and we have $\dim_k\H^1(G,\T_{\O})=t-1$.

\vspace{1ex}
(\ref{para-(..2)-2}.3)
If $p=2$, we use an ad-hoc construction. Let $\{ u_1 ,\dots , u_t \}$ be
a basis of $V$ and define a cocycle $\til{d}_0$ by $\til{d}_0(u_i):=u_i-u_i^2 x$ on basis elements and use
$\til{d}_0(u+v):=\til{d}_0u+(\til{d}_0v)^u$ inductively to define $\til{d}_0 u$ for all $u$. Notice that this requirement is compatible with conditions (i) ($\til{d}_0(2u)=0$) and (ii) (commutativity) from \ref{para-(..2)-1}. Thus we
get a well-defined cocycle on all of $V$ since $p=2$.  

The remaining part $\{a_0=a_1=0\}\cong\Hom_{\F_2}(V,k)$ contains
$2$-dimensional coboundaries spanned by $\iota$ and the Frobenius embedding
$$\mathrm{Frob}\colon V\hookrightarrow k \ : \  u\mapsto u^2.$$ Hence 
$$ \H^1(G,\T_{\O})\cong k\cdot[\til{d}_0] \ {\textstyle \bigoplus} \ 
\Hom_{\F_2}(V,k)/(k\cdot\iota+k\cdot\mathrm{Frob}),
$$
and, as the sum $k\cdot\iota+k\cdot\mathrm{Frob}$ is direct precisely when $t>1$, we get that $\dim_k\H^1(G,\T_{\O})$ is $t-1$ if $t>1$ or is $1$ of $t=1$.

\vspace{1ex}
\paragraph\label{para-(..3)}\ 
{\bf The case (\ref{para-cases}.3)}.\ 
Since $|H|$ is coprime to the characteristic of $k$, we have $$\H^1(G,M)=\H^1(N,M)^H  \leqno{\indent\textrm{{\rm (\ref{para-(..3)}.1)}}}
$$ (see, e.g., \cite{Br82}, III.10.4), 
The calculation of $\H^1(N,M)$ is similar to the one in the previous
paragraph, but $\F_p$ is replaced by $\F_q$ (for $q=p^s$) everywhere. The only difference
is for $p=2$, because then the subspace
$\Hom_{\F_{q}}(V,k)$ in $\Hom_{\F_p}(V,k)$ has a trivial 
intersection with the one dimensional subspace spanned by the Frobenius
embedding $u\mapsto u^2$.

We now describe the action of $H$ on 
the cohomology $\H^1(N,M)$.
First recall that the action of $H$ on $\T_{\O}$ is given by
$$(x^r\frac{d}{dx})^{\tau}=\zeta^{r-1}x^r\frac{d}{dx}.$$ 
It stabilizes $M$, on which is acts by 
$$ (a_0 + a_1 x + a_2 x^2)^\tau = \zeta^{-1} a_0 + a_1 x + \zeta a_2 x^2. $$ 
The action of $H$ on the cohomology 
$\H^1(N,M)$ is induced from the action on the space of cocycles given 
by $d^{\tau}u=(du^{\tau})^{\tau^{-1}}$ (cf.\ loc.\ cit.).
Hence if $p\neq 2$, we get $d^{\tau}_0u=(d_0\zeta u)^{\tau^{-1}}=\zeta^2d_0u+(\zeta-\zeta^2)ux-(\frac{1}{2}(\zeta-\zeta^2)u^2+\frac{1}{6}(1-\zeta^2)u)x^2$ where the last two terms form a coboundary; thus
$$
d^{\tau}_0=\zeta^2d_0,
$$
and if $p=2$, a similar result holds for the cocyle $\til{d}_0$ introduced in (\ref{para-(..2)-2}.3). This implies that the class $d_0$ is not $H$-invariant as long as $n\neq 2$, and $\til{d_0}$ is not $H$-invariant at all.
Next we look at the remaining part $\{a_0=a_1=0\}\cong\Hom_{\F_q}(V,k)$.
An element $a_2\in\Hom_{\F_q}(V,k)$ (corresponding to the cocycle
$du=a_2(u)x^2$) is $H$-invariant if and only if $d^{\tau}u=\zeta^{-1}
a_2(\zeta u)x^2=a_2(u)x^2$, or equivalently, $a_2(\zeta u)=\zeta a_2(u)$.
Since $\zeta$ generates $\F_{q}$, it is equivalent to the $\F_q$-linearity of
$a_2$.
Summing up, we have the following:

\vspace{1ex}
(\ref{para-(..3)}.1) Suppose $p\neq 2,3$. If $n\neq 2$, then 
$$
\H^1(G,\T_{\O})\cong\H^1(N,M)^H\cong\Hom_{\F_{q}}(V,k)/k\cdot\iota,
$$
and hence is of dimension $t/s-1$. If $n=2$ (which implies $s=1$), 
it is of dimension $t$, since $d_0$ is also $H$-invariant.

\vspace{1ex}
(\ref{para-(..3)}.2) If $p=3$, then we always have 
$$
\H^1(G,\T_{\O})\cong\H^1(N,M)^H\cong\Hom_{\F_{q}}(V,k)/k\cdot\iota,
$$
and $\dim_k\H^1(G,\T_{\O})=t/s-1$.

\vspace{1ex}
(\ref{para-(..3)}.3) If $p=2$, then $n\neq 2$ and $s>1$ (since $1 \neq n|2^s-1$).  Hence 
$$
\H^1(G,\T_{\O})\cong\H^1(N,M)^H\cong\Hom_{\F_{q}}(V,k)/k\cdot\iota,
$$
(recall that the Frobenius embedding is no longer present) and $\dim_k\H^1(G,\T_{\O})=t/s-1$.

\vspace{1ex}
\paragraph\label{infres}\ 
{\bf Remark.} \ Fix a subgroup $V_{1} = ({\Z}/p\Z)$ of $V$ and consider
the inflation restriction sequence for this subgroup (with $V':=V/V_1$):
$$ 0 \rightarrow H^1(V', M^{V_1}) \stackrel{\textrm{\footnotesize inf}}{\longrightarrow} H^1(V,M) \stackrel{\textrm{\footnotesize res}}{\longrightarrow} H^1(V_1,M)^{V'} $$
The left cohomology group is the following (for $p \neq 2$): the invariants $M^{V_1}$ are
just $kx^2$, on which $V'$ acts trivially, so it is isomorphic to
$\Hom(V',k)$, which maps via inflation to the part $\Hom_{\F_p}(V,k)/k\cdot\iota$ of 
the cohomology $H^1(G,\T_\O)$. On the other hand, the part $[d_0]$ is
precisely the one which is non-zero when mapped under restriction to 
$H^1(V_1,M)=H^1({\Z}/p\Z,\T_\O)$. This group is the one studied in \cite{BM00} (e.g. 4.2.2),
where an intrinsic characterization of
the class $[d_0]$ is provided. 

\vspace{2ex}
\sectioning\label{section-localmoduli}\ 
{\bf The local pro-representable hull}

\vspace{1ex}

In this section we will calculate 
the hull $H_{\rho}$ of $D_{\rho}$, where the group $G$ and its action on
$R=\O$ are as in (\ref{para-cases}.1), (\ref{para-cases}.2), and
(\ref{para-cases}.3).

The first one is a trivial case, i.e., the tame cyclic action is
rigid (compare \cite{BM00}, 4.3).

To analyze the case (\ref{para-cases}.2), we have to look at lifting
obstructions, much as is the case in \cite{BM00}, 4.2. In our case, liftings
are related to what we call formal truncated Chebyshev polynomials, which
we will now introduce first. 

\vspace{1ex}
\paragraph \label{chebyshev} {\bf Formal truncated Chebyshev polynomials.} \ 
For a positive integer $N$, we set
$$
M^{[N]}_{\alpha,\beta}(u)=\left[
\begin{array}{cc}
\sum^N_{k=0}{u+k-1\choose 2k}\alpha^k&
\alpha\sum^{N-1}_{k=0}{u+k\choose 2k+1}\alpha^k\\
\sum^{N-1}_{k=0}{u+k\choose 2k+1}\alpha^k+\beta(u)&
\sum^N_{k=0}{u+k\choose 2k}\alpha^k
\end{array}
\right],
$$
where $u$, $\alpha$ and $\beta(u)$ are indeterminates.
The entries are considered as formal power series of these indeterminates
with coefficients in $\Q$ (this means in particular that the binomial coefficients evaluate to polynomials in $u$).
We denote by $M^{[N]}_{\alpha}(u)$ $(1\leq N\leq \infty)$ the matrix
$M^{[N]}_{\alpha,\beta}(u)$ with $\beta(u)$ replaced by $0$.

As a matrix of formal power series, we have an 
identity
$$
M^{[\infty]}_{\alpha}(u)=\left[
\begin{array}{cc}
S_u(1+\frac{\alpha}{2})-(1+\alpha)S_{u-1}(1+\frac{\alpha}{2})&
\alpha S_{u-1}(1+\frac{\alpha}{2})\\
S_{u-1}(1+\frac{\alpha}{2})&
S_u(1+\frac{\alpha}{2})-S_{u-1}(1+\frac{\alpha}{2})
\end{array}
\right],
$$
where $S_u(x)$ is the Chebyshev polynomial of second kind
\begin{eqnarray*}
S_u(x)
&=&(u+1){}_2F_1(-u,u+2,\frac{3}{2};-\frac{x}{4}) \\
&=&\sin(u+1)\theta/\sin\theta\ \ (\textrm{if } u \in {\bf Z}_{\geq 0}, \textrm{ with } x = \cos \theta)
\end{eqnarray*}
Indeed, we have the following easy identity
$$
{\textstyle
S_u(1+\frac{x}{2})=\sum^{\infty}_{k=0}{u+k+1\choose 2k+1}x^k.
}
$$
For $1\leq N\leq\infty$ let us write
$$
M^{[N]}_{\alpha}(u)=\left[
\begin{array}{cc}
A^{[N]}_{\alpha}(u)&B^{[N]}_{\alpha}(u)\\
C^{[N]}_{\alpha}(u)&D^{[N]}_{\alpha}(u).
\end{array}
\right]
$$
Then it follows easily from the basic recursion between binomial coefficients that
$$
B^{[N]}_{\alpha}(u)=\alpha C^{[N]}_{\alpha}(u)\ \ \textrm{and}\ \
A^{[N]}_{\alpha}(u)+B^{[N]}_{\alpha}(u)=D^{[N]}_{\alpha}(u),
$$
i.e., the relations between the coefficients that hold for $N=\infty$ 
are also true for the truncated versions. 
It follows formally from these identites between entries that
$$
M^{[N]}_{\alpha}(u)M^{[N]}_{\alpha}(v)=M^{[N]}_{\alpha}(v)M^{[N]}_{\alpha}(u),
 \leqno{\indent\textrm{{\rm (\ref{chebyshev}.1)}}}
$$
as identities of matrices of formal power series in $u$, $v$ and 
$\alpha$. 

Next, we observe the following trigonometric identities:
\begin{condlist}
\item[(i)] $S_{u+v}(x)+S_{u-1}(x)S_{v-1}(x)=S_{u}(x)S_v(x)$,
\item[(ii)] $S_{u+v-1}(x)+2xS_{u-1}(x)S_{v-1}(x)=S_{u-1}(x)S_v(x)+S_u(x)S_{v-1}(x)$.
\item[(iii)] $S_u(x)^2-2xS_u(x)S_{u-1}(x)+S_{u-1}(x)^2=1$.
\end{condlist}
The first two of these imply that 
$$M^{[\infty]}_{\alpha}(u)M^{[\infty]}_{\alpha}(v)=M^{[\infty]}_{\alpha}(u+v).
 \leqno{\indent\textrm{{\rm (\ref{chebyshev}.2)}}}
$$
and the third one implies that 
$$
\det M^{[\infty]}_{\alpha}(u)=1.
 \leqno{\indent\textrm{{\rm (\ref{chebyshev}.3)}}}
$$

\vspace{1ex}
\paragraph \label{para-obstructions} {\bf Lifting obstructions.} \ Now assume $p \neq 2$, and let $A$ be an 
artinian local $k$ algebra with $A/\mathfrak{m}_A=k$.
Let $\alpha\in\mathfrak{m}_A$.
Let $V\subset k$ a finite dimensional $\F_p$-vector space, and 
$\beta\colon V\rightarrow\mathfrak{m}_A$ an $\F_p$-linear map.
For $u \in V \hookrightarrow k$ we use the following notation $$
\til{M}_{\alpha,\beta}(u)=M^{[\frac{p-1}{2}]}_{\alpha,\beta}(u)\ \ 
\textrm{and}\ \ \til{M}_{\alpha}(u)=M^{[\frac{p-1}{2}]}_{\alpha}(u)
$$
(where multiplication of $u$'s takes place inside $k$). By (\ref{chebyshev}.1) we have 
$$\til{M}_{\alpha}(u)\til{M}_{\alpha}(v)=
\til{M}_{\alpha}(v)\til{M}_{\alpha}(u)
\leqno{\indent\textrm{{\rm (\ref{para-obstructions}.1)}}}
$$
Also, by (\ref{chebyshev}.2) and (\ref{chebyshev}.3) we get
$$
\til{M}_{\alpha}(u)\til{M}_{\alpha}(v)\equiv\til{M}_{\alpha}(u+v)
\ \textrm{mod}\ (\alpha^{\frac{p-1}{2}})
\leqno{\indent\textrm{{\rm (\ref{para-obstructions}.2)}}}
$$
and 
$$
\det\til{M}_{\alpha}(u)\equiv 1\ \textrm{mod}\ (\alpha^{\frac{p-1}{2}+1}).
\leqno{\indent\textrm{{\rm (\ref{para-obstructions}.3)}}}
$$
We now look at what happens if $\beta \not \equiv 0$. A small calculation
shows that the commutation 
relation $
\til{M}_{\alpha,\beta}(u)\til{M}_{\alpha,\beta}(v)=
\til{M}_{\alpha,\beta}(v)\til{M}_{\alpha,\beta}(u)
$
is equivalent to 
$
\alpha\beta(u)C(v)=\alpha\beta(v)C(u).
$
Putting $v=1$ (where we tacitly assume to have conjugated the action of $V$ such that $\F_p \subseteq V$ (cf.\ \cite{BM00}, 4.2.1)) we get $\alpha\beta(u)=\alpha\beta(1)C(u)$. Thus, if $\alpha\beta$ is not a zero map, $C(u)$ must be a linear form, whence $\alpha=0$. So we get:

\vspace{1ex}
(\ref{para-obstructions}.4) \ For $\alpha\neq 0$, \ 
$\til{M}_{\alpha,\beta}(u)\til{M}_{\alpha,\beta}(v)=
\til{M}_{\alpha,\beta}(v)\til{M}_{\alpha,\beta}(u) \iff \alpha\beta=0$.
\vspace{1ex}

\noindent Looking at (\ref{para-obstructions}.3), we get

\vspace{1ex}
(\ref{para-obstructions}.5) \ If $\alpha \beta = 0$, then $\det\til{M}_{\alpha,\beta}(u)\equiv 1\ \textrm{mod}\ (\alpha^{\frac{p+1}{2}})$.
\vspace{1ex}

\noindent We can now prove the following:

\vspace{1ex}
{\bf \ref{para-obstructions}.6 Lemma.} \ {\sl The following conditions are equivalent:
\begin{condlist}
\item[(i)] $\til{M}_{\alpha,\beta}(u)\til{M}_{\alpha,\beta}(v)=
\gamma\til{M}_{\alpha,\beta}
(u+v)$ for all $u,v \in V$ and a $\gamma\in A$ (possibly depending on $u,v,\alpha,\beta\}$.
\item[(ii)] $\alpha^{\frac{p-1}{2}}=\alpha\beta=0$.
\end{condlist}
Actually, if this holds we have $\gamma=1$.}

\vspace{1ex}\noindent
\pf (i)$\Longrightarrow$(ii). The relation $\alpha\beta=0$ follows from (\ref{para-obstructions}.4). Since the identity in (i) obviously holds with $\gamma=1$ if $A=k$, 
we can set $\gamma=1+\delta$ with $\delta\in\mathfrak{m}_A$.
Taking determinants, it follows from (\ref{para-obstructions}.3) that 
$\gamma^2(=1+\delta(2+\delta))=1+P\alpha^{\frac{p+1}{2}}$ for some $P \in A$.
Since $\delta+2$ is invertible in $A$, 
$$
\gamma\equiv 1 \mbox{ mod }\alpha^{N+1}
\leqno{\indent\textrm{{\rm (\ref{para-obstructions}.7)}}}
$$ 
for $N=\frac{p-1}{2}$. 
The coefficient of $\alpha^{\frac{p-1}{2}}$ in the lower-left entry of 
$\til{M}_{\alpha,\beta}(u)\til{M}_{\alpha,\beta}(v)$ is
$$
{\textstyle
\sum^{N-1}_{k=0}{u+k\choose 2k+1}{v+N-k-1\choose 2(N-k)}+
\sum^{N-1}_{k=0}{v+k\choose 2k+1}{u+N-k\choose 2(N-k)},}
\leqno{\indent\textrm{{\rm (\ref{para-obstructions}.8)}}}
$$
(note that we have used that we know already that $\alpha\beta=0$ so there is no contribution from
$\beta$ in this calculation.)
We now put $u=N$ and $v=2$ (again assuming tacitly that $\F_p \subseteq V$) in (\ref{para-obstructions}.8). The result is $1$, so the above coefficient (\ref{para-obstructions}.8) is non-zero.
Since the coefficient of $\alpha^{\frac{p-1}{2}}$ in the lower-left entry of $\til{M}_{\alpha,\beta}(u+v)$ is zero,  we conclude  from this and (\ref{para-obstructions}.7) that $\alpha^N=0$
as desired.

(ii)$\Longrightarrow$(i). If $\beta \equiv 0$ this is already in (\ref{para-obstructions}.2), so let $\beta \not \equiv 0$. Using $\alpha\beta=0$ and $A(u)\equiv 1$ mod $\alpha$, we calculate the entries of the matrix $M:=\til{M}_{\alpha,\beta}(u)\til{M}_{\alpha,\beta}(v)$ as 
follows:
\begin{eqnarray*} M_{1,1}
&=& A(u)A(v)+\alpha C(u)C(v) \\
 M_{1,2}
&=& \alpha A(u)C(v)+\alpha A(v)C(u)+\alpha^2 C(u)C(v) \\
 M_{2,1}
&=& A(v)C(u)+A(u)C(v)+\alpha C(u)C(v)+\beta(u+v) \\
 M_{2,2} 
&=& A(u)A(v)+\alpha A(u)C(v)+\alpha A(v)C(u)+\alpha(1+\alpha)C(u)C(v) 
\end{eqnarray*}
Note that, except for the third one, these entries are independent of $\beta$.
Since we already have (\ref{para-obstructions}.2), we find that these
entries are equal to 
$A(u+v)$, $\alpha C(u+v)$, $C(u+v)+\beta(u+v)$ and $A(u+v)+\alpha C(u+v)$,
respectively.
Thus we have the relation in (i) with $\gamma=1$.
\qed
 
\vspace{1ex}
\paragraph \label{liftings} {\bf Explicit liftings.} \ Assume $p \neq 2$. We try to lift the action of $V$
on $\T_\O$ to $A$ using the frational linear transformation corresponding to 
the truncated formal Chebyshev polynomials:
$$ x^u := \frac{ax+b}{cx+d}, \mbox{ where } {{a\,b}\choose{c\,d}} = \til{M}_{\alpha,\beta}(-u)
\leqno{\indent\textrm{{\rm (\ref{liftings}.1)}}}$$
for $\beta: V \rightarrow \mathfrak{m}$ a linear map and $\alpha \in A$ (notice
the sign change from $u$ to $-u$, since we are lifting $x^u=x/(1-ux)$). 

At the infinitesimal level $A=k[\epsilon]$ with $\epsilon^2=0$, we let $\alpha=a_0 \epsilon$ and
$\beta(u)= - \epsilon \cdot \phi(u)$ for a linear map $\phi \in \Hom_{\F_p}(V,k)$. Lemma \ref{para-obstructions}.6 assures us of the
fact that (\ref{liftings}.1) does define a first order lifting as long as
$\alpha=0$ for $p=3$. It is explicitly given as
$$ x^u = \frac{ (1+\frac{1}{2}u(u+1)a_0 \epsilon)x-ua_0 \epsilon } { 1-(u+\frac{1}{6}u(u^2-1)a_0 \epsilon - \epsilon \phi(u))x+\frac{1}{2}u(u-1)a_0 \epsilon} $$ 
The corresponding cocycle $d$ is by the formula in proposition 
\ref{pro-lifting} given as
\begin{eqnarray*}
du &=&
\frac{d}{d\epsilon}\left(
\frac{x^u}{1+ux^u}\right)\bigg|_{\epsilon=0}  \\ &=& -a_0 u + a_0(u^2+u) x -(a_0(\frac{1}{3}u^3+\frac{1}{2}u^2+\frac{1}{6}u)- \phi(u)) x^2 \\
&=& a_0 d_0 + \phi 
\end{eqnarray*}
in $Z^1(V,M) = k d_0 + \Hom_{\F_p}(V,k)$, where $d_0$ is as in (\ref{para-(..2)-2}.1). This means that (\ref{liftings}.1) for $\alpha=0$ defines a lifting
in the direction of $\phi$ (which is unobstructed) and for $\beta=0$ in the direction of $a_0 [d_0]$ (obstructed by $\alpha^{\frac{p-1}{2}}=0$). If both
$\alpha$ and $\beta$ are non-zero, a lifting in the direction of $a_0[d_0]+\phi$
is obstructed by the equations in \ref{para-obstructions}.6.

\vspace{1ex}
\paragraph \label{hull} {\bf Calculating the hull.} \ 

\vspace{1ex}
{\bf \ref{hull}.1 The case $p \neq 2,3$}. \ If we let $R$ be
the ring
$$ R = k[[x_0,x_1,\dots ,x_t]]/ \langle x_0^{\frac{p-1}{2}},x_0 x_1, x_0 x_2, \dots , x_0 x_{t}, x_1+ \dots +x_t \rangle, $$ 
then \ref{para-obstructions}.6 shows that
(\ref{liftings}.1) defines a lifting of $\rho$ to $R$, and hence there
is a morphism of functors 
$$ \Hom(R, -) \rightarrow D_\rho. $$
To prove that $R$ is actually the hull $H_\rho$ of $D_\rho$,  we argue
as in \cite{BM00}, p.\ 217.
It 
suffices to prove that $R$ is a versal deformation, i.e., that the above morphism is smooth and that it is an isomorphism on the level of tangent spaces. The latter is clear from our computation of the tangent space to $D_\rho$ from 
group cohomology and the above explicit form of $R$. To prove the former, let $A' \rightarrow A$ be a small 
extension in $\C_k$ with kernel $I$. 
We have to show that 
$$ \Hom(R,A') \rightarrow D(A') \times_{D(A)} \Hom(R,A) $$
is surjective.
So assume that $\rho^\sim \in \mbox{im}(\Hom(R,A) \rightarrow D(A))$ lifts
to $\rho^\approx$ in $D(A')$. This means the corresponding obstruction
in $\H^2(G,\T_\O) \otimes I$ is zero. Since $G=({\bf Z}/p{\bf Z})^t$, this obstruction measures exactly the possible failure of the commutation relation $\rho^\approx(u+v)
= \rho^\approx(u) \rho^\approx(v), \forall u,v \in V$, which we know by
(\ref{para-obstructions}.6) is given by the equations in $R$. So we can
lift via a $\til{M}_{\alpha',\beta'}(u)$ to $\Hom(R,A')$, and we 
can adjust this lifting in such a way that its image in $D(A')$ coincides with
the given one, since the tangent spaces to the two functors are isomorphic (note that the fibers are $H^1(G,\T_\O) \otimes I$-torsors).
This finishes the proof of smoothness.

\vspace{1ex}
{\bf \ref{hull}.2 The case $p = 3$}. \ The same argument works, except 
that the class $[d_0]$ does not occur so that all liftings in the direction
of $\Hom_{\F_p}(V,k)$ are unobstructed. The result is
$$ H_\rho = k[[x_1, \dots ,x_{t}]]/\langle x_1+ \dots +x_t \rangle. $$

\vspace{1ex}
{\bf \ref{hull}.3 The case $p = 2$}. \ We will deal with this case by an
ad-hoc construction. Recall that in this case $G \cong V \subseteq k$ is a $t$-dimensional $\F_2$-vector space for which we pick a basis $\{ u_1, \dots, u_t \}$. As before, let $\alpha \in A$ and $\beta: V \rightarrow \mathfrak{m}_A$ be a linear map. We lift
the action of this {\sl basis} $u_i$ by 
$$ x^{u_i} := \frac{x+\alpha u_i}{(u_i+\beta(u_i))x+1}. $$
We have $(x^{u_i})^{u_i}=x$, and we observe that $(x^{u_i})^{u_j}=(x^{u_j})^{u_i}$ is equivalent
to 
$$\alpha u_i \beta(u_j) = \alpha u_j \beta(u_i)
\leqno{\indent\textrm{{\rm (\ref{hull}.4)}}}$$
for all $i,j$. We now {\sl define} a lift to any element of $V$ by 
$$ x^{u_{i_1}+\dots+u_{i_l}} := x^{u_{i_1} \circ \dots \circ u_{i_l}} $$
for $i_1,\dots,i_l \in \{1,\dots,t\}$. This forces commutativity to hold.
Notice that this only gives a well-defined lift to all of $V$ since $p=2$. 
If we set $A=K[\epsilon],\alpha=a_0 \epsilon,\beta=\epsilon \phi(u)$ for a
linear map $\phi \in \Hom_{\F_p}(V,k)$, then we find that this lift
corresponds to the cocycle $a_0 d_0 + \phi \in Z^1(V,M)$, so that we
are indeed lifting in all directions of the tangent space. 
One can now reason as before to find that the hull of $D_\rho$ is 
given by 
$$ H_\rho = k[[x_0,x_1,\dots,x_t]]/ \langle  x_1+\cdots+x_t,u_1x_1+
\cdots+u_tx_t , x_0 (x_i u_j - x_j u_i)_ {i,j=1 \dots t}  \rangle $$

{\footnotesize
\unitlength0.7pt
\begin{center}
\begin{picture}(300,150)
\thicklines
\multiput(0,50)(40,20){2}{\line(0,1){50}}
\multiput(0,50)(0,50){2}{\line(2,1){40}}
\put(20,85){\vector(1,0){30}}
\put(57,82){$[d_0]$-direction}
\put(-50,125){$\Hom(V,k)$-}
\put(-50,110){directions}
\put(5,30){$p \neq 2$}

\multiput(200,50)(40,20){2}{\line(0,1){50}}
\multiput(200,50)(0,50){2}{\line(2,1){40}}
\put(205,75){\line(1,0){30}}
\put(235,95){\line(1,0){40}}
\put(235,75){\line(2,1){40}}
\put(268,78){$[d_0]$-direction}
\put(150,125){$\Hom(V,k)$-}
\put(150,110){directions}
\put(205,30){$p = 2$}

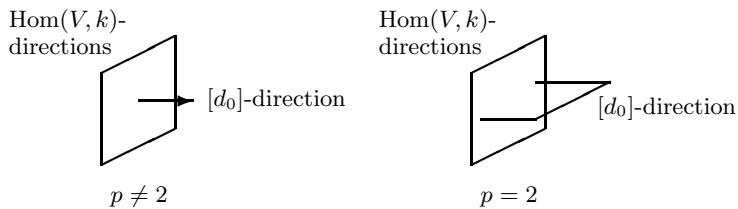
\put(-35,0){Pictorial representation of the local versal deformation ring}

\end{picture}
\end{center}}

{\bf \ref{hull}.5 The case $n \neq 1$.} \  
Finally, in the general case (\ref{para-cases}.3),
the arguments are completely similar to those of the preceeding paragraphs;
the results are as follows:

\vspace{1ex}
(\ref{hull}.5.1) If $n\neq 2$ or $p=2,3$, then only first-order deformations
comming from $\Hom_{\F_q}(V,k)$ occur, and these can be lifted without 
obstruction: 
$$
H_{\rho}\cong k[[x_1,\ldots,x_{t/s}]]/\langle x_1+\cdots+x_{t/s} \rangle.
$$

({\ref{hull}.5.2) If $p \neq 2,3$ and $n=2$ (hence $s=1$), then the action of $H$ extends to the lifting
(\ref{liftings}.1) (just replacing $u$ by $\zeta u$), and the 
obstructions do not change. Hence we have
$$
H_{\rho}\cong k[[x_0,x_1,\ldots,x_t]]/\langle x_0^{\frac{p-1}{2}},x_0 x_1, \dots, x_0 x_t, x_1+\cdots+x_t \rangle.
$$

\vspace{1ex}
\paragraph\label{local-theorem}\ {\bf Theorem.} \ {\sl Let $\rho : G \rightarrow \mbox{Aut}(\T_\O)$ be a local representation of a finite group $G$, where $\O$ is
of characteristic $p$. Let $n$ be an integer coprime to $p$, define
$s:= \min \{ s'  :  n | p^{s'}-1 \}$, and let $[d_0]$ be the
cohomology class defined in \textup{(\ref{para-(..2)-2})}. The following table 
lists the dimension of the group cohomology $H^1(G,\T_\O)$, the fact
whether $[d_0]$ is trivial \textup{(---)}, unobstructed \textup{(unobs.)} or leads to obstructions \textup{(obs.)}, 
and the Krull-dimension $\dim H_\rho$ of the pro-representable hull $H_\rho$ of the local
deformation functor $D_\rho$: }

\begin{center}
\begin{tabular}{ll|l|c|l}
$G$ & $(p,t,n)$ & $h^1(G,\T_\O)$ & $[d_0]$ & $\dim H_\rho$ \\
\hline 
$\Z/n$ & & 0 & --- & 0 \\
$(\Z/p)^t$ & $p \neq 2,3$ & $t$ & obs. & $t-1$ \\
& $p=3$ & $t-1$ & --- & $t-1$ \\
& $p=2,t>1$ & $t-1$ & obs. & $t-2$ \\
& $p=2,t=1$ & 1 & unobs. & 1 \\
$(\Z/p)^t \sd \Z/n$  & $n \neq 2$ or $p=2,3$ & $t/s-1$ & --- & ${t/s-1}$ \\
& $n=2$ & $t$ & obs. & $t-1$ \\
\end{tabular}
\end{center} 

\vspace{1ex}
{\bf \ref{local-theorem}.1 Remark.} \ For $n=t=1$, this result agrees with proposition 4.1.1 in \cite{BM00}, where it is shown that
$$ h^1(G,\T_\O) = \lfloor \frac{2 \beta}{p}  \rfloor - \lceil \frac{\beta}{p} \rceil $$
for $G=\Z/ p$ a cyclic $p$-group and $\beta = \sum_{j=0}^\infty (|G_i|-1)$
(recall that $G_i$ are the higher ramification groups). Indeed, in our case $G_0=G_1=G$
and $G_i=0$ for $i>1$, so that $\beta=2p-2$ and we get
$$ h^1(G,\T_\O) = \left\{ \begin{array}{ll} 3-2=1 & \mbox{ if } p>3, \\ 2-2 = 0 & \mbox{ if } p=3, \\  2-1 = 1  & \mbox{ if } p=2. \end{array} \right. $$
Similarly, our calculation of the hulls and its Krull-dimension is compatible
with the results from \cite{BM00} if we observe that (in their notation, cf.\ p.\ 215) $\psi(X) = X^{\frac{p-1}{2}}$ mod $p$ ($p>2$). 

\vspace{1ex}
{\bf \ref{local-theorem}.2 Remark.} \ Obstructions to lifting first order
deformations are certain second cohomology classes in $\H^2(G,\T_\O)$, 
but can form a strict subset of this cohomology group. If $t=1$, then
the group is easy to calculate directly (cf.\ \cite{BM00}), or seen 
to be isomorphic to $\H^1(G,\T_\O)$ by Herbrand's theorem since $G$ is cyclic.
This already shows that for $t=1$, obstructions can form only a small
part of $\H^2(G,\T_\O)$. The computation of this second cohomology group
in the general case $t>1$ (so $G$ is no longer cyclic) seems rather tedious.

\vspace{2ex}

\sectioning\label{section-main-algebraic}\ 
{\bf Main theorem on algebraic equivariant deformation}

\vspace{1ex}

We have now collected all information needed to prove our main algebraic theorem:

\vspace{1ex}
\paragraph\label{global-theorem}\ {\bf Theorem.} \ {\sl Let $X$ be 
an ordinary curve over a field of characteristic $p>0$ and let
$G$ be a finite group of acting on $X$ via $\rho  :  G \rightarrow \mbox{Aut}(X)$. 

\textup{(a)} The Krull-dimension
of the pro-representable hull of the deformation functor $D_{X,\rho}$ 
is given by
$$ \dim H_{X,\rho} = 3 g_Y-3 +\delta + \sum_{i=1}^s \dim H_{\rho_i}, $$
where $g_Y$ is the genus of $Y:=G \backslash X$, $\delta$ is given in \textup{(\ref{para-branchpoints})} and
$H_{\rho_i}$ is the pro-representable hull of the local
deformation functor associated to the representation of the ramification
group $G_i$ at the branch point $w_i$ in $X \rightarrow Y$, whose dimension was given in \textup{(\ref{local-theorem})}, \textup{unless} in the following four
cases:
\begin{enumerate}
\item when $p=2, Y={\bf P}^1$ and $X \rightarrow Y$ is branched above
2 points; then $\dim H_{X,\rho} = \dim H_{\rho_1} + \dim H_{\rho_2}$; 
\item when $X = {\bf P}^1 \rightarrow Y = {\bf P}^1$ is tamely branched above two points, then $\dim H_{X,\rho}=0$;
\item   when $X = {\bf P}^1 \rightarrow Y = {\bf P}^1$ is wildly branched above a unique point; then  $\dim H_{X,\rho} = \dim H_{\rho_1}$, (in this case, 
$G$ is a pure-$p$ group and $X \rightarrow Y$ is an Artin-Schreier cover); 
\item when  $X \rightarrow Y$ is an unramified cover of elliptic
curves; then $\dim H_{X,\rho} =1$.
\end{enumerate}
Furthermore, if the genus $g$ of $X$ is $\geq 2$, then $D_{X,\rho}$ is
pro-representable by $H_{X,\rho}$. 

\medskip

\textup{(b)} The dimension of the tangent space to the functor
$D_{X,\rho}$ as a $k$-vector space satisfies
$$ \dim_k D_{X,\rho}(k[\epsilon]) = \dim H_{X,\rho} + \left\{ \begin{array}{ll} \# \{ i : n_i \leq 2 \} & \mbox{ if } p \neq 2,3; \\ 0 & \mbox{ if } p=3; \\ \# \{ i : n_i=1 \mbox{ and } t_i>1 \} & \mbox{ if } p=2. \end{array} \right. $$
}

\vspace{1ex}
\pf 
Let $Y$ be the quotient $G \backslash X$.
From proposition (\ref{pro-globalhull}), we find that 
$$ \dim H_{X,\rho} = \sum_{i=1}^s \dim H_{\rho_i} + h^1(Y,\pi^G_\ast \T_X), $$
where $H_{\rho_i}$ are the hulls of the local deformation functors
associated to the action of the ramification groups $G_i$ at
wild ramifications points $w_1,...,w_s$ on the space of local
derivations, which was computed in (\ref{local-theorem}). 

To compute the $h^1$-term, recall from (\ref{pro-linebundle}) that $\pi^G_\ast \T_X=\T_Y(-\Delta)$, where $\Delta$ is defined in (\ref{para-branchpoints}). By Riemann-Roch, we find
$$ h^1(\pi^G_\ast \T_X) = 3g_Y-3+\delta + h^0(\T_Y(-\Delta)), $$
where $g_Y$ is the genus of $Y$. 

Since $\deg(\T_Y(-\Delta))=2-2g_Y-\delta$, the last term vanishes if $g_Y>1$
or $g_Y=1$ and $\delta>0$ or $g_Y=0$ and $\delta>2$.

If $g_Y=1$ and $\delta=0$, $X$ is an unramified cover of an elliptic
curve, hence is an elliptic curve itself. 

Assume that $g_Y=0$, that there are at least two branch points on $Y$
and that $\delta \leq 2$.
If $p \neq 2$, then these branch points have to be tame, so $\delta=2$, and
the Hurwitz formula implies that $g_X=0$ too.
If $p=2$, they can be wild, but both ramification groups have to be $\Z/2\Z$, so we still have $\delta=2$; in both cases, $h^0(\T_Y(-\delta)) =h^0(\mathcal{O}_{{\bf P}^1})=1$.

If $g_Y=0$ and only one point on $Y$ is branched, then it follows
from Hurwitz's formula (using the fact that second ramification groups
vanish in the ordinary case, cf.\ \cite{Nakajima:87}) that $g_X=0$ too, and
the ramification has to be wild at this point. So if $p \neq 2$ or $p=2, t>1$, 
we find $\delta + h^0(\T_Y(-\delta)) = 2 + h^0(\mathcal{O}_{{\bf P}^1})=3$.
On the other hand, if $p=2$ and $t=1$, then $\delta + h^0(\T_Y(-\delta))  =1 + h^0(\mathcal{O}_{{\bf P}^1}(1))=3$. 
Let $np^t$ be the order of the ramification group at that unique
point, where $n$ is coprime to $p$. Hurwitz's formula gives in particular
that $(n-1)p^t+2$ divides $2np^t$, and this (together with $n|p^t-1$) 
implies $n=1$. Hence we do get a ${\bf Z}/p{\bf Z}$-cover.
This finishes the proof of part (a). 

For part (b), we apply the formula from (\ref{globaltangent}) in combination
with (\ref{para-spectral}.1). It thus suffices to compute $h^0(Y, \mathcal{H}^1(G,\T_X))$, but $\mathcal{H}^1(G,\T_X)$ is concentrated in the 
branch points $w_i$, where it equals the group cohomology $H^1(G_i,\T_{\O_{w_i}})$ (\cite{BM00}, 3.3), so $\dim H_{X,\rho}$ and
$\dim_k D_{X,\rho}(k[\epsilon])$ differ only at places where 
$[d_0]$ is obstructed in table (\ref{local-theorem}). \qed

\vspace{1ex}
{\bf \ref{global-theorem}.1 Remark.} \ The main algebraic theorem stated
in the introduction follows from \ref{global-theorem} by excluding the cases
$g \leq 2, p = 2, 3$. 

\vspace{1ex}
{\bf \paragraph\label{def-example}.1 Example (Artin-Schreier curves).} \  The Artin-Schreier
curve whose affine equation is given by $(y^{p^t}-y)(x^{p^t}-x)=c$ for
some constant $c \in k^*$ has automorphism group $$G=(\Z/p\Z)^{2t} \sd
D_{p^t-1},$$ where $D_\ast$ denotes a dihedral group of order $2 \ast$. The quotient $Y=G \backslash
X$ is a projective line, and the branching groups are $\Z/2\Z$ (twice if $p \neq 2$ and once if $p=2$) and $(\Z/p\Z)^t \sd \Z/(p^t-1)\Z$ (once). This
curve with its full automorphism groups hence allows for a one-dimensional
deformation space (for different reasons if $p\neq 2$ and $p=2$). This 
deformation is exactly given by varying $c$.

\vspace{1ex}
{\bf \ref{def-example}.2 Example (Drinfeld modular curves).} \ The Drinfeld 
modular curves $X(n)$ from the introduction have automorphism
group $G:=\Gamma(1)/\Gamma(n)\F_q^*$ for $d:=\deg(n)>1$. The quotient $Y:=G \backslash X(n)$ is a projective line, over which $X$ is branched at 2 points with 
ramification groups $\Z/(p+1)\Z$ and $\F_q^d \sd \F_q^*$ respectively. 
Hence $X(n)$ can be deformed in $d-1$ ways (regardless of $p$, but
for different reasons if $p=2$ -- then we are in exceptional 
case (1) from (\ref{global-theorem})).  

\vspace{3ex}

{\bf PART B: ANALYTIC THEORY}

\vspace{2ex}

\sectioning\label{section-teich}\ 
{\bf Equivariant deformation of Mumford curves}

\vspace{1ex}
\paragraph\label{para-mumford} {\bf Mumford curves.} \ Let $(K,|\cdot|)$ be a complete discrete valuation field with valuation ring 
$\O_K$ and residue field $\O_K/\m_{\O_K}=k$. Recall that a projective curve $X$ over $K$ is called a {\sl Mumford curve} if
it is ``uniformized over $K$ by a Schottky group''. This means that there exists
a free subgroup $\Gamma$ in $PGL(2,K)$ of rank $g$, acting on ${\bf P}^1_K$ with limit set ${\cal L}_\Gamma$ such that $X$ satisfies $X^\an \cong \Gamma
\backslash ({\bf P}^{1,\an}_K - {\cal L}_\Gamma)$ as rigid analytic
spaces. Mumford (\cite{Mumford:72}) has shown that these conditions
are equivalent to the existence of a stable model of $X$ over $\O_K$
whose special fiber consists only of rational components with $k$-rational double points. Because of the ``GAGA''-correspondence for one-dimensional rigid analytic spaces, we do not have to (and will not) distinguish between analytic and algebraic curves. 
It is well-known that 
Mumford curves are ordinary (this is basically because
their Jacobian is uniformized by $({\bf G}^{\an}_{m,K})^g/\Gamma^{\mbox{{\rm {\tiny ab}}}}$ where $g$ is the genus of $X$, cf.\ \cite{Cornelissen:01}, 1.2). Thus, the results from 
the previous section essentially solve the equivariant deformation 
problem for Mumford curves in a cohomological way. In this part however, 
we want to develop an independent theory of analytic deformation of 
Mumford curves based on the groups that uniformize them. This will 
make the liftings and obstructions whose cohomological existence was
proven in the previous part more ``visible'' as actual deformations
of $2 \times 2$-matrices over $K$.

\vspace{1ex}
 \paragraph\label{para-mumauto} {\bf Automorphisms.} \ It is well-known
(\cite{Cornelissen:01}, 1.3) that for a Mumford curve $X$ of genus $g
\geq 2$ with Schottky group $\Gamma$, 
$\A(X) = N(\Gamma)/\Gamma$, where $N(\Gamma)$ is the normalizer of $\Gamma$
in $PGL(2,K)$. 
Conversely, if $N$ is a discrete subgroup of $PGL(2,K)$ containing $\Gamma$ and contained in $N(\Gamma)$, then it induces a group 
of automorphisms $$\rho: N/\Gamma \hookrightarrow \A(X).$$

\vspace{1ex}
 \paragraph\label{para-nothom*} {\bf Notation.} \ 
If a finitely generated discrete subgroup $N$ of $PGL(2,K)$ is given,
let $\Hom^*(N,PGL(2,K))$ denote the set of {\sl injective}
homomorphisms $\phi : N \rightarrow PGL(2,K)$ {\sl with discrete image}. 
Then such $N$ contains a finite index normal free subgroup of finite
rank $\Gamma$ (\cite{Gerritzen:80}, I.3),
and if $\Gamma$ is non-trivial, the pair $(N,\Gamma)$ gives rise to a
Mumford curve with an action of $N/\Gamma$ as is being considered here.

\vspace{1ex}
 \paragraph\label{para-mumrigid} {\bf Rigidity.} \ Two Mumford curves $X,X'$
with Schottky groups $\Gamma, \Gamma'$ are isomorphic if and only if 
$\Gamma$ and $\Gamma'$ are conjugate in $PGL(2,K)$ (\cite{Mumford:72}, 4.11). 

\vspace{1ex}
{\bf \ref{para-mumrigid}.1 Remark.} \
Note that this is very different from the situation in the uniformization 
theory of Riemann surfaces $S$, where in a representation $S=\Gamma \backslash
\Omega$ with $\Gamma$  a Schottky subgroup of  $PGL(2,{\bf C})$, the domain 
of discontinuity $\Omega$ of $\Gamma$ is not the universal topological covering space
of $S$, whereas this does hold for Mumford curves. 

\vspace{1ex}
\paragraph\label{para-tangenttoM} {\bf Analytic
deformation functors.} \ 
Recall that $\C_K$ is the category of local Artinian $K$-algebras.
If $N$ and $\phi \in \Hom^*(N,PGL(2,K))$ are given, we consider the 
{\sl analytic deformation functor} $$D_{N,\phi}: \C_K \rightarrow
{\tt Set} $$ of the pair $(N,\phi)$, which sends $A \in \C_K$ to
the set
of liftings of $(N,\phi)$ to $A$. Here, a lifting is a morphism
$\phi^\sim \in \Hom(N, PGL(2,A))$ which, when composed with reduction
modulo the maximal ideal $\m_A$ of $A$, equals the original embedding
$\phi$ (in particular, $\phi^\sim$ is injective). Note that we do not consider classes of liftings modulo
conjugacy by $PGL(2)$ --- this implies that $D_{N,\phi}$ is naturally
equipped with an action of group functor $PGL(2)^\wedge$ given
by 
$$ PGL(2)^\wedge : \C_K \rightarrow \mbox{\tt{Groups}} : A \mapsto \ker [PGL(2,A) \rightarrow PGL(2,K)], $$  and we denote the quotient by $$D_{N,\phi}^\sim
:= PGL(2)^\wedge \backslash D_{N,\phi}. $$

Since $N$ is finitely generated, it is not difficult to show that the
functors $D_{N,\phi}$ and $D^\sim_{N,\phi}$ 
have pro-representable hulls $H_{N,\phi}$ and $H^\sim_{N,\phi}$. We want to compute the
dimension
of the tangent spaces to these functors, and the Krull-dimension of
their hulls, and this will be done by ``decomposing'' $N$ using its
structure
as a group acting on a tree to give a decomposition of $D_{N,\phi}$; note that 
such a decomposition is {\sl not} given on the level of $D^\sim_{N,\phi}$.

\vspace{2ex}

\sectioning\label{section-structureBT}\ {\bf Structure of $N$ as a group acting on a tree.} \ 

\vspace{1ex}

We fix
a finitely generated discrete subgroup $N$ of $PGL(2,K)$, and we will now
recall how the structure of $N$ can be seen from its action on the Bruhat-Tits
tree (cf.\ section 2 of \cite{Cornelissen:01}). 

\vspace{1ex}
\paragraph\label{para-structureBT1} {\bf The Bruhat-Tits tree.} \
Let
$\T$ denote the Bruhat-Tits tree of $PGL(2,K)$ ({\sl i.e.,} its vertices are similarity classes $\Lambda$ of rank two ${\cal O}_K$-lattices in $K^2$, and two vertices are connected  by an edge if the corresponding quotient module has length one -- {\sl see} Serre \cite{Serre:80}, Gerritzen \& van der Put \cite{Gerritzen:80}). We assume $K$ to be large enough so that all fixed
points of $N$ are defined over $K$; then $N$ acts without inversion on
$\T$. It is a regular
tree in which the edges emanating from a given vertex are in one-to-one
correspondence with ${\bf P}^1(k)$. The tree $\T$ admits a left action
by $PGL(2,K)$.

\vspace{1ex}
\paragraph\label{para-structureBT2} {\bf Notations on trees.} \
For any subtree $T$ of $\T$, let $\Ends(T)$ denote its set of ends ({\sl i.e.,}
equivalence classes of half-lines differing by a finite
segment). There is a natural correspondence between ${\bf P}^1(K)$ and
$\Ends(\T)$.  Let $V(T)$ and $E(T)$ denote the set of vertices and
edges of $T$ respectively.  For $\sigma \in E(T)$, let $o(\sigma)$ (respectively
$t(\sigma)$) denote the origin (respectively, terminal) vertex of
$\sigma$.  Let $N_x$ denote the stabilizer of a vertex or edge $x$ of $T$ for the action of $N$. The maps $o,t$ induces maps $N_\sigma \rightarrow N_\Lambda$ for $\Lambda=o(\sigma)$ or $\Lambda = t(\sigma)$, which will be denoted by the same letter. For any $u,v \in {\bf P}^1(K)$, let $]u,v[$ denote
the apartment in $\T$ connecting $u$ and $v$ (seen as ends of $\T$). 

\vspace{1ex}
\paragraph\label{para-structureBT4} {\bf The trees associated to $N$.} \
We can construct a locally finite tree  $\T({\cal L})$ (possibly empty) from any compact subset ${\cal L}$ of ${\bf P}^1(K)$: it is the minimal subtree of $\T$ whose set of ends coincides with ${\cal L}$, or equivalently, the minimal subtree of $\T$ containing  
$ \bigcup_{u,v \in {\cal L}} \ ]u,v[. $

We define $\T_N$  to be the tree associated to the subset ${\cal L}_N$
consisting of the limit points of $N$ in ${\bf P}^1(K)$. Since $N$ is a finitely generated discrete group, $\T_N$ coincides with the tree of $N$ as it is defined in Gerritzen \& van der Put \cite{Gerritzen:80}.

\vspace{1ex}
\paragraph\label{para-structureBT5} {\bf The graph associated to $N$.} \
$\T_N$ admits a natural action of $N$, and we denote the quotient graph by
$T_N:= N \backslash \T_N$; the corresponding quotient map will be denoted by $\pi_N$. The graph $T_N$ is finite and connected.

We turn $T_N$  into a graph of groups as follows:
let $T$ be a spanning tree (maximal subtree) of $T_N$, which we can
see  as a
subtree of $\T_N$ by a fixed section $\iota: T \rightarrow \T_N$ of $\pi_N$. Let $c=c(T_N)$ denote
the cyclomatic number of $T_N$ (= number of edges outside
$T$), and fix $2c$ lifts $e_i^\pm$ of these edges outside $T$ to
$\T_N$ which satisfy: $t(e_i^+) \in V(\iota (T)), o(e_i^-) \in V(\iota
(T))$. Fix $c$ hyperbolic elements $\{ \gamma_i \}_{i=1}^c$ in $N$
such that $\gamma_i e_i^+ = e_i^-$. Then $\iota(T) \cup \{ e_i^+
\}_{i=1}^c $ is a fundamental domain for the action of $N$ on $\T_N$. 

For  any vertex  $\Lambda
\in V(\T_N)$ and edge  $\sigma =[\Lambda,M] \in E(\T_N)$ we denote 
by $N_\Lambda$ and $N_\sigma =N_\Lambda \cap N_M$ their respective
stabilizers for the action of $N$. Note that these groups are finite since $N$ is discrete.

For a vertex $v \in V(T_N)=V(T)$, we let $N_v = N_{\iota(v)}$. For
edges $e \in E(T_N)$, either $e \in E(T)$, and then we let
$N_e=N_{\iota(e)}$, or else, there is a unique $i$ such that
$\pi_N(e^\pm_i)=e$, and we let $N_e=N_{e_i^+}$. 

The morphisms between these groups are defined as follows: if $e \in
E(T)$, then $N_e \hookrightarrow N_{t(e)}$ and $N_e \hookrightarrow N_{o(e)}$ are the natural
inclusions; if, on the other hand, $e=\pi_N(e^\pm_i)$, then $N_e
\hookrightarrow N_{t(e)}$ is the natural inclusion, but $N_e
\hookrightarrow N_{o(e)}$ is given by $s \mapsto \gamma_i^{-1} s \gamma_i$. 

We then
have the following description of the group $N$:

\vspace{1ex}
\paragraph\label{para-structureBT6} {\bf Theorem.} (Bass-Serre
\cite{Dicks:89}, 4.1, 4.4, \cite{Serre:80}) {\sl For any spanning tree $T$ of $T_N$,  $N$ equals the fundamental
group of the graph of groups $T_N$ at $T$. This means that $N$ is generated 
by the amalgam of $N_v$ over $N_e$ for all $e \in E(T_N), v \in V(T_N)$ together
with the fundamental group of $T_N$ at $T$ as a plain graph, viz., the
free group $F_c$ on $c$ generators $\{ n_i \}_{i=1}^c$, where $c=c(T_N)$ is
the cyclomatic number of $T_N$. The further relations  in $N$ are of the form $n_i t(\gamma) n_i^{-1} = o(\gamma)$
for every $i=1,...,c$ and for every $\gamma \in N_e, e \in T_N -
T$. In particular, there is a split exact sequence of groups
$$ 0 \rightarrow \lim_{ \stackrel{\longrightarrow}{T_N^\sim} }
N^\sim_\bullet \rightarrow N {\rightarrow} F_{c(T_N)} \rightarrow 0, $$
where $\pi : T^\sim_N \rightarrow T_N$ is the universal covering of $T_N$ as a plain graph,
which has been made into a graph of groups by setting $N^\sim_\bullet$ for $\bullet \in V(T^\sim_N) \cup E(T^\sim_N)$
equal to $N_{\pi(\bullet)}$. \qed} 

\vspace{2ex}

\sectioning\label{section-finitecase}\ {\bf Decomposition of the
  functor $D_{N,\phi}$.} 

\vspace{1ex}
\paragraph\label{prop-tofinite} {\bf Proposition.} \ {\sl Let $s:
  F_{c(T_N)} \rightarrow N$ be a splitting of the sequence in
  (\ref{para-structureBT6}). Then there is an isomorphism of functors
$$ D_{N,\phi} \cong \lim_{\stackrel{\longleftarrow}{T_N}}
D_{N_\bullet,\phi|_{N_\bullet}} \times D_{F_c,\phi \circ s}, $$ where
the inverse limit is in the category of functors (note that morphisms
between $N_\bullet$ naturally induce morphisms of functors between $D_{N_\bullet}$). 
}

\vspace{1ex}
{\bf \ref{prop-tofinite}.1 Remark.} \ Note that we get a direct
product of functors, but a limit of functors over $T_N$ (instead of
the obvious semi-direct
product and limit over $T_N^\sim$). We also note that there is no such decomposition on the level of the functors $D^\sim_{N,\phi}$. 

\vspace{1ex}
\pf  Let $A \in \C_K$.  By restriction, a deformation of $N$ to $A$ trivially gives rise to deformations of $N_\bullet$
and $F_c$.  

For the rest of the proof, we will imitate the construction of $T_N$
as a graph of groups, but we will lift to $T^\sim_N$ instead of
$\T_N$. So choose a fixed maximal spanning tree $ \iota : T
\hookrightarrow T^\sim_N$ and a basis $\{ \gamma_1,...,\gamma_c \}$ of
$s(F_c)$, where $c=c(T_N)$. Take, as before, $2c$
edges $e_i^\pm \in E(T^\sim_N)$ such that  $t(e^+_i) \in V(T), o(e^-_i)
\in V(T), \gamma_i e^+_i = e^-_i$. Thus, $T \cup \{ e_i^+ \}$ is a
``fundamental domain'' for $ T^\sim_N \rightarrow T_N$. 

To give elements in 
$\lim\limits_{\stackrel{\longleftarrow}{T_N}}
D_{N_\bullet,\phi|_{N_\bullet}}(A) \mbox{ and } D_{F_c,\phi \circ s}(A)$
means precisely to give a compatible collection of $ \phi_v : N_v \hookrightarrow PGL(2,A)$
and $\phi_c: F_c \hookrightarrow PGL(2,A)$. Compatibility means that for
$e \in E(T_N),$ the following diagram is commutative:
$$ \begin{array}{ccl} N_e & \longrightarrow & N_{o(e)} \\
                       \bigdownarrow & & \bigdownarrow
                       \rlap{$\scriptstyle{\phi_{o(e)}}$} \\
N_{t(e)} & \underrel{\longrightarrow}{\phi_{t(e)}} & PGL(2,A). 
\end{array} $$
We want to extend this to an
embedding of $N$.

By the construction of the fundamental domain, there exists for any $v \in V(T^\sim_N)$, a unique $\gamma \in
s(F_c)$ such that $v \in \gamma T$, and this allows us to define $N_v
\rightarrow PGL(2,A)$ to be $ \sigma \mapsto \phi_{\gamma^{-1} v} (\gamma^{-1} \sigma \gamma)$.
For edges $e$, we similarly get $\gamma$ such that $e \in \gamma \cdot
(T \cup \{ e_i^+ \}_{i=1}^c)$, and the same works. 

By compatibility, we thus get an embedding of
$\lim\limits_{\stackrel{\longrightarrow}{T_N^\sim}} N_\bullet$ into
  $PGL(2,A)$, which by construction is compatible with the conjugation
  action of $F_c$, so that we finally get an embedding $N \hookrightarrow
  PGL(2,A)$, viz., an element of $D_{N,\phi}(A)$.  Since this
  construction is functorial in $A$, we get the desired inverse map of
  functors. \qed

\vspace{1ex}
\paragraph\label{para-freecase} {\bf Computing the functor
  $D_{F_c,\phi \circ s}$.} \ The set of morphisms 
$$ \Hom(F_c, PGL(2,K)) $$
is a smooth algebraic variety over $K$; by choosing a basis of $F_c$, it is
isomorphic 
to $PGL(2,K)^c$ over $K$. We can take its formal 
completion at the $K$-rational point $\phi \circ s$, and
$$ D_{F_c,\phi \circ s} \cong \Hom(F_c, PGL(2,K))^\wedge_{\phi \circ s} $$
as formal functors. In particular, $$ \dim_K D_{F_c,\phi \circ s}(K[\epsilon])  =
\dim H_{F_c, \phi \circ s} =  3c, $$
where the first one is the dimension of the tangent space and the
second one the Krull-dimension of the pro-representable hull of the functor.

\vspace{1ex}
\paragraph\label{para-finitecase} {\bf Computing the functor
  $D_{N,\phi}$ for finite $N$.} \ 
Here, the argument is based on the simple observation that an injective
element $\phi$ of $\Hom(N,PGL(2,K))$ corresponds to a cover ${\bf P}^1 \rightarrow {\bf P}^1$ with Galois group $N$; and hence, it is related to the algebraic
deformation functor (\ref{global-def}) of the pair $({\bf P}^1,\phi)$ (regarding $\phi$ as a 
representation of $N$ into $\mathrm{Aut}({\bf P}^1)$). The functor
$D_{{\bf P}^1, \phi}$ is defined modulo conjugation by $PGL(2)$,
whereas the analytic deformation functor $D_{N,\phi}$ carries a
natural action of $PGL(2)^\wedge$. However, it is easy to see that
$$ D^\sim_{N,\phi} = D_{{\bf P}^1, \phi}. $$
From this formula, we get in particular that 
$$
\dim H_{N,\phi} =  \dim H_{{\bf P}^1,\phi} + 3 - \nu(\phi(N)),
\leqno{\indent\textrm{{\rm (\ref{para-finitecase}.1)}}}
$$
where for a finite subgroup $G \subseteq PGL(2,K)$,
$$ \nu(G) = \dim \mathrm{Nor}_{PGL(2,K)}(G) $$
is the dimension of the normalizer of $G$ in $PGL(2,K)$ as an
algebraic group. Formula (\ref{para-finitecase}.1) continues to hold when
hulls are replaced by tangent spaces.

\vspace{1ex}
By Dickson's {\sl Hauptsatz}, the finite groups $N$ acting on ${\bf
  P}^1$ in positive characteristic are
known. Let us first set up the notation:

\vspace{1ex}
{\bf \ref{para-finitecase}.2 Notation.} \ We let $D_n$ denote the dihedral group of order $2n$. We will write $P(2,q)$ to denote either $PGL(2,q)$ or $PSL(2,q)$ by slight abuse of notation, with the convention that any related numerical 
quantities that appear between set delimiters $\{ \}$ are only to be considered for $PSL(2,q)$. 

\vspace{1ex}

We now recall this classification in the version as it is given in Valentini-Madan, as this more geometrical form immediately allows
us to compute $D_{N,\phi}$ using the results from section \ref{section-main-algebraic}:

\vspace{1ex}
{\bf \ref{para-finitecase}.3 Theorem} (Dickson, cf.\ \cite{Madan:80}). \ {\sl Any finite subgroup of $PGL(2,K)$ is isomorphic to a finite subgroup of $PGL(2,p^m)$ for some $m>0$. The group $PGL(2,p^m)$ has the following finite subgroups $G$, such that $\pi_G$ is branched over $d$ points with ramification groups isomorphic to
$G_1,...,G_d$:

(i) $G={\bf Z}/n\Z$ for $(n;p)=1$, $d=2$, $G_1=G_2={\bf Z}/n\Z$;

(ii) $G=D_n$ with $p \neq 2$, $n|p^{m} \pm 1$, $d=3$, $G_1=G_2={\bf Z}/2\Z,G_3={\bf Z}/n\Z$ or also, $p=2$, $(n;2)=1$, $d=2$ and $G_1={\bf Z}/2\Z ,G_2={\bf Z}/n\Z$; 

(iii) $G=(\Z/p\Z)^t \sd \Z/n\Z$ for $t \leq m$ and
$n|p^m-1, n|p^t-1$ with $d=2$ and $G_1=G,G_2={\bf Z}/n\Z$ if $n>1$ and $d=1,G_1=G$ otherwise;

(iv) $G=P(2,p^t)$ with $d=2$ and $G_1=(\Z/p\Z)^t \sd \Z/\{\frac{1}{2}\}(p^t-1)\Z,G_2={\bf Z}/{\{\frac{1}{2}\}(p^t+1)}\Z$; 

(v) $A_4$ of $p \neq 2,3$, $d=3$, $G_1={\bf Z}/2\Z, G_2=G_3={\bf Z}/3\Z$;

(vi) $S_4$ if $p \neq 2,3$, $d=3$, $G_1={\bf Z}/2\Z, G_2=G_3={\bf Z}/4\Z$;

(vii) $A_5$ if $5|p^{2m}-1$ and $p \neq 2,3,5$ with $d=3$ and $G_1={\bf Z}/2\Z,G_2={\bf Z}/3\Z,G_3={\bf Z}/5\Z$ or $p=3$, $d=2$ and $G_1=\Z/3\Z \sd \Z/2\Z,G_2={\bf Z}/5\Z$. \qed}

\vspace{1ex}
{\bf \ref{para-finitecase}.4 Lemma.} \ {\sl The normalizer of a finite
subgroup $N$ of $PGL(2,K)$ has dimension $\nu(N)=0$, unless if $N$ is 
cyclic of order prime-to-$p$; then $\nu(N)=1$, or if $N$ is a pure
$p$-group; then $\nu(N)=2$.}

\vspace{1ex}
\pf Any group from the above list which does not belong to the
mentioned exceptions has at least three fixed points on ${\bf P}^1$,
the set of which should also remain stable under the action of the normalizer of
$N$, which hence is finite.

A cyclic subgroup $N$ of order prime-to-$p$ has a diagonalizable
generator, and by a direct computation, this is seen to be exactly stabilized by the one-dimensional group $D$ generated by the
center of $PGL(2,K)$ and the involution ${0\,1}\choose{1\,0}$.

A $p$-group $N$ can be put into upper diagonal form by conjugation,
and a little computation shows that the 
stabilizer of such a group consists precisely of the 2-dimensional group of upper trigonal matrices.  \qed

\vspace{1ex}
\paragraph \label{main-analytic} {\bf Theorem.} \ {\sl Let $N$ be as in section \ref{section-structureBT}, and suppose a Bass-Serre representation of $N$ 
is given as in \ref{para-structureBT6}. Then
$$ \dim  H^\sim_{N,\phi} = 3 c(T_N)-3 + \sum_{v \in V(T_N)} h(N_v) - \sum_{e \in  E(T_N)} h(N_e), $$
and
$$ \dim_K D^\sim_{N,\phi}(K[\epsilon]) = 3 c(T_N)-3 + \sum_{v \in V(T_N)} t(N_v) - \sum_{e \in  E(T_N)} t(N_e), $$
where for a finite group $G \subset PGL(2,K)$, the numbers $h(G)$ and
$t(G)$ are given in the table below:}

$$ \begin{tabular}{ll|l|l} 
$G$ & $(p,t,n)$ & $h(G)$ & $t(G)$ \\
\hline
$\Z/n\Z$ & $(n;p)=1$ & $2$ & $2$ \\
$D_n$ & $p\neq 2$ & $3$ & $3$ \\
        & $p=2$ & $4$ & $4$ \\
$(\Z/p\Z)^t$ & $p \neq 2,3$ & $t$ & $t+1$ \\
             & $p=3$  & $t$ & $t$ \\
             & $p=2, t>1$ &  $t-1$ & $t$ \\
             & $p=2, t=1$ & $2$ & $2$ \\
$(\Z/p\Z)^t \sd \Z/n\Z$ & $p\neq 2 \ \& \ n\neq 2$ or $p=2,3$ & $t/s+2$ & $t/s+2$ \\
       & $n=2$ & $t+2$ & $t+3$ \\
$P(2,p^t)$ & $\{ p^t \neq 5 \}$ & $3$ & $3$ \\
$A_4,S_4$ & & $3$ & $3$ \\
$A_5$ & $p \neq 3$ & $3$ & $3$ \\
    & $p = 3$ & $3$ & $4$ 
\end{tabular} $$
\pf   Since $N$ is infinite, the action of $PGL(2)^\wedge$ is of
dimension 3. By \ref{prop-tofinite} and \ref{para-freecase}, we are
reduced to computing $D_{N,\phi}$ for finite $N$ occuring in $T_N$. We
know which different $G$ can occur on the edges and vertices of $T_N$
by Dickson's theorem \ref{para-finitecase}.3. For each of these, using
\ref{para-finitecase}.1 we are reduced to the computation of the
algebraic data, for which we appeal to \ref{global-theorem}, and to
the computation of $\nu(G)$, which is in lemma
\ref{para-finitecase}.4. If we let $h^{\mathrm{alg}}(G)$ and
$t^{\mathrm{alg}}(G)$ denote the Krull-dimension of the
pro-representable hull and vector space dimension of the tangent space to $D_{{\bf P}^1,\phi|G}$ respectively, we find
$$ \begin{tabular}{ll|l|l|l} 
$G$ & $(p,t,n)$ & $h^{\mathrm{alg}}(G)$ & $t^{\mathrm{alg}}(G)$ & $3-\nu(G)$\\
\hline
$\Z/n\Z$ & $(n;p)=1$ & $0$ & $0$ & $2$ \\
$D_n$ & $p\neq 2$ & $0$ & $0$ & $3$ \\
        & $p=2$ & $1$ & $1$ & $3$ \\
$(\Z/p\Z)^t$ & $p \neq 2,3$ & $t-1$ & $t$& $1$ \\
             & $p=3$  & $t-1$ & $t-1$ & $1$ \\
             & $p=2, t>1$ &  $t-2$ & $t-1$ & $1$ \\
             & $p=2, t=1$ & $1$ & $1$ & $1$ \\
$(\Z/p\Z)^t \sd \Z/n\Z$ & $p\neq 2 \ \& \ n\neq 2$ or $p=2,3$ & $t/s-1$ & $t/s-1$ & $3$\\
       & $n=2$ & $t-1$ & $t$ & $3$ \\
$P(2,p^t)$ & $\{ p^t \neq 5 \}^\ast$ & $0$ & $0$ & $3$ \\
$A_4,S_4$ & &  $0$ & $0$ & $3$ \\
$A_5$ & $p \neq 3$ & $0$ &  $0$ & $3$ \\
    & $p = 3$ & $0$ & $1$ & $3$ 
\end{tabular} $$
In this computation, note at $^\ast$ that $PSL(2,5) = A_5$ does not occur in Dickson's list if $p=5$. \qed

\vspace{1ex}
{\bf \ref{main-analytic}.1 Remark.} \ The main analytic theorem stated
in the introduction follows from \ref{main-analytic} by excluding the cases
$g \leq 2, p = 2, 3$. 

\vspace{2ex}

\sectioning\label{section-vs}\ 
{\bf Compatibility between algebraic and analytic deformation}

\vspace{1ex}
\paragraph\label{para-analyticdef}\ {\bf Deformation of Mumford uniformization.} \ 
We have already seen how to compare analytic and algebraic deformation functors for finite groups acting on ${\bf P}^1$. We now want to compare these functors in general, in particular to achieve equality between the apparently different results from the main algebraic and analytic theorem (\ref{global-theorem} and \ref{main-analytic}). 

 Let $X$ be a Mumford curve over $K$ uniformized 
by a Schottky group $\Gamma$ and $ \phi : N \hookrightarrow PGL(2,k)$ a discrete group between $\Gamma$
and its normalizer in $PGL(2,K)$. Let $\rho : N/\Gamma \hookrightarrow \A(X)$. To be able to compare the functors $D^\sim_{N,\phi}$ and
$D_{X,\rho}$, it will be necessary to develop a theory of 
algebraic deformation of Mumford uniformization; this means just to 
translate the formalism of Mumford (\cite{Mumford:72}) from fields $K$ to 
elements in the category $\C_K$. For lack of a reference, we sketch it here; the reader who wants to follow the details in encouraged to take a copy of \cite{Mumford:72} at hand. 

\vspace{1ex}
\paragraph\label{para-Artinian}\ {\bf Analytic objects in $\C_K$.} \
The maximal ideal of an object $A$ in $\C_K$ is denoted by $\m_A$.
Each object $A$ in $\C_K$ can be made into a $K$-affinoid algebra in a unique way by a 
suitable surjective homomorphism $K\langle X_1,\ldots,X_n\rangle
\rightarrow A$ over $K$ (cf.\ \cite{BGR} (6.1)).
We denote by $\O_A$ (resp.\ $\m_{\O_A}$) the subring (resp.\ the ideal in
$\O_A$) consisting of power-bounded (resp.\ topologically nilpotent) 
elements in $A$ (cf.\ loc.\ cit.\ (6.2.3)).
Since every element in $\m_A$ is nilpotent, we have $\O_A\bigcap\m_A\subseteq
\m_{\O_A}$. By this, it is easily seen that $\O_A$ is a local ring with
the maximal ideal $\m_{\O_A}$, that $\O_A/\m_{\O_A}\cong k$, and that
$\pi^{-1}_A(\O_K)=\O_A$ and $\pi^{-1}_A(\m_{\O_K})=\m_{\O_A}$, where 
$\pi_A\colon A\rightarrow K$ is the reduction map.

\vspace{1ex}

{\bf Example.} \ In the ring of dual numbers $K[\epsilon]$, the ring of 
power-bounded elements is $\O_K+K\epsilon$, whereas the ideal of topologically
nilpotent elements is $\m_{\O_K}+K\epsilon$.

\vspace{1ex}
\paragraph\label{para-lattice}\ {\bf Lattices.} \ 
Let $A$ be an object in $\C_K$. By a lattice in $A^2$ we mean an 
$\O_A$-submodule $M$ in $A^2$ that is free of rank $2$. By  elementary commutative
algebra, this is equivalent to $M\subset A^2$ being an $\O_A$-submodule 
such that the image $\ovl{M}$ in $K^2$ by the reduction map $A^2\rightarrow
K^2$ is a lattice in the usual sense.
We consider the set $\Delta^{(0)}_A$ of similarity classes of lattices
up to multiplication by $A^{\ast}$. Then $\Delta^{(0)}_A$ can 
be naturally identified with the
set of equivalence classes of couples $(\P,\phi)$, where $\P$ is an $\O_A$-scheme isomorphic to $\P^1_{\O_A}$ and $\phi$ is an isomorphism between
$\P\otimes A$ and $\P^1_A$, and two couples $(\P,\phi)$ and
$(\P',\phi')$ are equivalent if there exists an $\O_A$-isomorphism
$\psi\colon\P\rightarrow\P'$ such that $\phi'\circ\psi=\phi$. 
The identification between $\Delta^{(0)}_A$ and the 
space of such couples is given by $$M\mapsto\P(M)=\Proj(\Sym_{\O_A}M)$$ where
$\phi$ is induced from $M\otimes A\cong A^2$.

\vspace{1ex}
\paragraph\label{para-flat}\ {\bf Trees.} \ 
We take a subgroup $N\subset\PGL(2,A)$ such that its image $\ovl{N}$ in
$\PGL(2,K)$ by the reduction map is finitely generated, discrete and isomorphic to $N$.
Such a subgroup $N$ contains a normal free subgroup $\Gamma$ of finite index, 
since $\ovl{N}$ does and $N$ and $\ovl{N}$ are isomorphic as groups. 
This $\Gamma$ satisfies a ``flatness'' condition analogous to 
\cite[(1.4)]{Mumford:72} (or, equivalently, property $\ast$ in 
loc.\ cit.\ pp.139) in the folowing sense:
if $\Sigma$ is the set of all sections $\Spec \ A\rightarrow\P^1_A$ fixed by
non-trivial elements $\gamma\in\Gamma$, then 
for any $P_1,P_2,P_3,P_4\in\Sigma$, the cross-ratio $R:=R(P_1,P_2;P_3,P_4)$ or 
its inverse $R^{-1}$ lie in $\O_A$
(note: $\Sigma$ does not depend on the choice of $\Gamma$ in $N$).
The proof is easy from the fact $\pi^{-1}_A(\O_K)=\O_A$.
Given $P_1,P_2,P_3$ with homogeneous coordinates $w_1,w_2,w_3$, respectively, 
let $M=\O_Aa_1w_1+\O_Aa_2w_2+\O_Aa_3w_3$, where the $a_i$ satisfy a non-trivial 
linear relation $a_1w_1+a_2w_2+a_3w_3=0$.
The class $v(P_1,P_2,P_3)$ of $M$ in $\Delta^{(0)}_A$ depends only on 
$P_1,P_2,P_3$.
We let $\Delta^{(0)}_{\Gamma}$ be the
set of all such $v(P_1,P_2,P_3)$.
The set $\Delta^{(0)}_{\Gamma}$ is ``linked'' in the sense of 
loc.\ cit.\ (1.11), and the thus obtained tree is obviously the usual tree 
with respect to $\ovl{\Gamma}$.

\vspace{1ex}
\paragraph\label{para-join}\ 
The construction of the formal scheme also parallels the original one.
For $M_1$ and $M_2$ in $\Delta^{(0)}_{\Gamma}$ one defines the join 
$\P(M_1)\vee\P(M_2)$ to be the closure of the graph of the birational map
$\P(M_1)\cdots\rightarrow\P(M_2)$ induced from $\phi_2^{-1}\circ\phi_1$,
where $(\P(M_i),\phi_i)$ corresponds to $M_i$ ($i=1,2$) by the
correspondence from (\ref{para-lattice}).
The formal scheme ${\cal P}_{\Gamma}$ over $\Spf\,\O_A$ is then 
constructed as in loc.\ cit.\ pp.156 using these joins.
Obviously, its fiber over 
$\Spf\,\O_K$ is isomorphic to the usual formal scheme; in particular, their
underlying topological spaces are isomorphic.
It is clear that the associated rigid space $\Omega_{\Gamma}$ of 
${\cal P}_{\Gamma}$ in the sense of \cite{BL93} \S 5 is the complement 
in $\P^{1,\mathrm{an}}_A$ of the closure of the set of fixed rig-points
(corresponding to the fixed sections).
The quotient and the algebraization can equally well be taken by a r 
reasoning similar to the usual case. What finally comes out is a scheme 
$X_{\Gamma}$ over $A$ with special fiber over $K$ the Mumford curve
corresponding to $\ovl{\Gamma}$, and hence one can further take a
finite quotient by $N/\Gamma$. 

\vspace{1ex}
\paragraph\label{para-functors}\ 
The above construction of infinitesimal deformation of Mumford uniformization
induces a morphism of functors 
$$
\Phi \ \colon \ 
D^\sim_{N,\phi}\longrightarrow
D_{X,\rho},
$$
by associating to a deformation of $N \hookrightarrow PGL(2,K)$ to
$\bar{N} \hookrightarrow PGL(2,A)$ the corresponding ``Mumford'' curve over $\Spec\, A$. 
By an argument parallel to \cite{Mumford:72} \S 4, it is not difficult
to see the following:

\vspace{1ex}
\paragraph\label{prop-functors} {\bf Proposition.} \ 
{\sl The morphism $\Phi$ is an isomorphism.}
\qed

\vspace{1ex}
\paragraph\label{prop-functcor} {\bf Remark.} \ If $X$ is a
  Mumford curve uniformized by a Schottky group $\Gamma$ and $N$ is
  between $\Gamma$ and its normalizer in $PGL(2,K)$, let $\rho \colon
  N\backslash \Gamma \rightarrow \Aut(X)$ the corresponding
  representation. Then the (a priori very different looking)
  results from the algebraic computation
  \ref{global-theorem} for $D_{X,\rho}$ and the analytic computation
  \ref{main-analytic} for $D^\sim_{N,\phi}$ agree. There might be a
  more direct combinatorial proof of this equality.

\vspace{1ex}
{\bf \ref{example2}.1 Example (Artin-Schreier-Mumford curves).} \ If, for the Artin-Schreier curves of
\ref{def-example}.1 over a non-arch\-ime\-de\-an field $K$ the value of $c$ satisfies $|c|<1$, then $X_{t,c}$
is a Mumford curve, and the corresponding normalizer of its Schottky
group is 
$$ N_t = (({\bf Z}/p\Z)^t \sd {\bf Z}/(p^t-1)\Z) \ast_{{\bf
    Z}/{(p^t-1)}\Z} D_{p^t-1}.$$  
We compute from this that
the analytic infinitesimal deformation space is $1$-dimen\-si\-o\-nal, in
concordance with the algebraic result.

\vspace{1ex}
\paragraph\label{example2}{\bf .2} {\bf Example (Drinfeld modular curves).} \ The Drinfeld modular curve $X({\mathfrak n})$ is known
to be a Mumford curve (cf.\ \cite{Gekeler:96}), and the normalizer of its Schottky group is
isomorphic to an amalgam (cf.\ \cite{Cornelissen:01}) $$ N({\mathfrak n}) =  PGL(2,p^t) \ast_{({\bf Z}/p\Z)^t \sd {\bf Z}/{(p^t-1)\Z}} ({\bf Z}/p\Z)^{td} \sd {\bf Z}/{(p^{t}-1)}\Z $$
(at least if $p \neq 2, q \neq 3$). The above formula gives a
$(d-1)$-dimensional infinitesimal analytic deformation space, and this agrees with the algebraic result. 

\vspace{2ex}

\begin{small}

\noindent {\bf Acknowledgments.} The first author is post-doctoral fellow of the Fund for Scientific Research - Flanders (FWO -Vlaanderen). This work was done while the first author was visiting Kyoto University and the MPIM in Bonn, and while the second author was visiting the University of Paris-VI. 
Thanks to Frans Oort for pointing us to \cite{BM00}.

\vv

\vv

\noindent Max-Planck-Institut f{\"u}r Mathematik, Vivatsgasse 7,
53111 Bonn, Germany (gc {\sl current)}

\vv

\noindent University of Gent, Dept.\ of Pure Mathematics, 
Galglaan 2, 9000 Gent, Belgium (gc)

\vv

\noindent Kyoto University, Faculty of Science, Dept.\ of Mathematics, Kyoto 606-8502, Japan
\vv

\noindent e-mail: {\tt gc@cage.rug.ac.be}, {\tt kato@kusm.kyoto-u.ac.jp}

\end{small}

\end{document}